\documentclass[11pt]{article}
\usepackage{amssymb}
\usepackage{amsmath}
\usepackage{latexsym}

\input{prepictex}
\input{pictex}
\input{postpictex}
\begin{document}
\bibliographystyle{plain}
\baselineskip=14pt
\newtheorem{lemma}{Lemma}[section]
\newtheorem{theorem}[lemma]{Theorem}
\newtheorem{prop}[lemma]{Proposition}
\newtheorem{cor}[lemma]{Corollary}
\newtheorem{definition}[lemma]{Definition}
\newtheorem{definitions}[lemma]{Definitions}
\newtheorem{remark}[lemma]{Remark}
\newtheorem{conjecture}[lemma]{Conjecture}
\newcommand{\TT}[1]{\widetilde{T}\mbox{}_{#1}}
\newcommand{\bb}{{\bf B}'({\bf d})\cup {\bf B}''({\bf d})}
\newcommand{\vanish}{{\bf B}\setminus \Psi(\bb)}
\newcommand{\word}[1]{s_{i_1}s_{i_2}\cdots s_{i_{#1}}}
\newcommand{\dd}{\mbox{\scriptsize\bf d}}
\newcommand{\E}[1]{E_{i_1}E_{i_2}\cdotsE_{i_{#1}}}
\newcommand{\F}[1]{F_{i_1}F_{i_2}\cdotsF_{i_{#1}}}
\newcommand{\EA}[1]{E_{i_1}^{(r_1)}E_{i_2}^{(r_2)}\cdots
           E_{i_{#1}}^{(r_{#1})}}
\newcommand{\FA}[1]{F_{i_1}^{(r_1)}F_{i_2}^{(r_2)}\cdots
           F_{i_{#1}}^{(r_{#1})}}
\newcommand{\UA}{U_{\cal A}^-}
\newcommand{\dimr}{\left( \begin{array}{c} n+1 \\ r \end{array} \right)}
\newcommand{\TE}{\widetilde{E}}
\newcommand{\TF}{\widetilde{F}}
\newcommand{\C}{{\mathbb{C}}}
\newcommand{\N}{{\mathbb{N}}}
\newcommand{\Q}{{\mathbb{Q}}}
\newcommand{\R}{{\mathbb{R}}}
\newcommand{\Z}{{\mathbb{Z}}}
\newcommand{\B}{{\mathbf{B}}}
\newcommand{\CA}{{\mathcal{A}}}
\newcommand{\CL}{{\mathcal{L}}}
\newcommand{\0}{-\,}
\newcommand{\cdi}{\mbox{CD$(\mathbf{i})$}}
\newcommand{\cdii}{\mbox{CD$(\mathbf{i}')$}}

\begin{center}
\baselineskip=13pt
{\LARGE Rectangle diagrams for the Lusztig cones of quantized enveloping
algebras of type $A$}

\vspace*{10mm}
{\Large Robert Marsh}

{\em Department of Mathematics and Computer Science, University of Leicester,
University Road, Leicester LE1 7RH, England} \\
{E-mail: R.Marsh@mcs.le.ac.uk}

\baselineskip=21pt

{\bf Abstract} \\
\parbox[t]{5.15in}{
Let $U$ be the quantum group associated to a Lie algebra
$\bf g$ of type $A_n$. The negative part $U^-$ of $U$ has
a canonical basis $\mathbf{B}$ defined by Lusztig and Kashiwara, with
favourable properties.
We show how the spanning vectors of the cones defined by
Lusztig~\cite{lusztig7}, when regarded as monomials in Kashiwara's root
operators, can be described using a remarkable rectangle combinatorics.
We use this to calculate the Lusztig parameters of the corresponding
canonical basis elements, conjecturing
that translates of these vectors span the simplicial regions of linearity of
Lusztig's piecewise-linear function~\cite[\S2]{lusztig2}. \\
\parskip=14pt
\noindent {\bf Keywords:}
quantum group, canonical basis, Lusztig cone, longest word, piecewise-linear
combinatorics.
} \\ \ \\
\end{center}


\section{Introduction}

Let $U=U_q({\mathbf{g}})$ be the quantum group associated to a semisimple
Lie algebra $\mathbf{g}$ of rank $n$. The negative part $U^-$ of $U$ has
a canonical basis $\mathbf{B}$ with favourable properties
(see Kashiwara~\cite{kash2} and Lusztig~\cite[\S14.4.6]{lusztig6}).
For example, via action on highest weight vectors it gives rise to bases
for all the finite-dimensional irreducible highest weight $U$-modules.

Let $W$ be the Weyl group of ${\mathbf{g}}$, with Coxeter generators
$s_1,s_2,\ldots ,s_n$, and let $w_0$ be the element of maximal length in $W$.
Let ${\mathbf{i}}$ be a reduced expression for $w_0$, i.e. $w_0=s_{i_1}s_{i_2}
\cdots s_{i_k}$ is reduced. Lusztig obtains a parametrization of the
canonical basis $\mathbf{B}$ for each such reduced expression $\mathbf{i}$,
via a correspondence between a basis of PBW-type associated to $\mathbf{i}$
and the canonical basis. This gives a bijection
$\varphi_{\mathbf{i}}:\B\rightarrow \N^k,$ where $\N=\{0,1,2,\ldots \}$.

Kashiwara, in his approach to the canonical basis (which he calls the global
crystal basis), defines certain root operators $\TF_i$ on the canonical
basis (see~\cite[\S3.5]{kash2}) which lead to a parametrization of the
canonical basis for each reduced expression $\mathbf{i}$ by a certain subset
$Y_{\mathbf{i}}$ of $\N^k$. This gives a bijection $\psi_{\mathbf{i}}:\B\rightarrow
Y_{\mathbf{i}}$. The subset $Y_{\mathbf{i}}$ is called the string cone.

There is no simple way to express the elements of $\B$ in terms of the
natural generators $F_1,F_2,\ldots ,F_n$ of $U^-$. This has been done for
all of $\mathbf{B}$ only in types $A_1,A_2,A_3$ and $B_2$
(see~\cite[\S3.4]{lusztig2},~\cite[\S13]{lusztig7},~\cite{xi2}
and~\cite{xi3}) and appears to become arbitrarily complicated in general.

A monomial $F_{i_1}^{(a_1)}F_{i_2}^{(a_2)}\cdots F_{i_k}^{(a_k)}$, where
$F_i^{(a)}=F_i^a/[a]!$, is said to be {\em tight} if it belongs to ${\B}$.
Lusztig~\cite{lusztig7} described a method which in low rank cases leads
to the construction of tight monomials. He defined, for each reduced
expression $\mathbf{i}$ of $w_0$, a certain cone $C_{\mathbf{i}}$ in $\N^k$
which we shall call the Lusztig cone.
Let $$M_{\mathbf{i}}=\{F_{i_1}^{(a_1)}F_{i_2}^{(a_2)}\cdots F_{i_k}^{(a_k)}
\,:\,\mathbf{a}\in C_{\mathbf{i}}\}$$
be the set of monomials obtained from elements of $C_{\mathbf{i}}$.
Lusztig showed that, in types $A_1,A_2$ and $A_3$,
$M_{\mathbf{i}}\subseteq \B$. The author~\cite{me7} has extended this to type
$A_4$. However, counter-examples of Xi~\cite{xi4} and Reineke~\cite{reineke2}
show that this fails for higher rank.

Instead, one can consider monomials in Kashiwara's operators $\TF_i$. Define,
for a reduced expression $\mathbf{i}$ for $w_0$,
$$\B_{\mathbf{i}}=\{b\in \B\,:\,b\equiv \TF_{i_1}^{a_1}\TF_{i_2}^{a_2}\cdots
\TF_{i_k}^{a_k}\cdot 1,\ \mathbf{a}\in C_{\mathbf{i}}\}.$$
It has been shown by Premat~\cite{premat1} that
$C_{\mathbf{i}}\subseteq Y_{\mathbf{i}}$, so we have
$\B_{\mathbf{i}}=\psi_{\mathbf{i}}^{-1}(C_{\mathbf{i}})$.
Then it is known that $\B_{\mathbf{i}}=M_{\mathbf{i}}$ in types $A_1,A_2,A_3$
(using Lusztig's work~\cite{lusztig7} and direct calculation); it will be
seen in~\cite{me10} that this is also true in type $A_4$.
This breaks down in higher ranks because of the above
counter-examples. Thus $B_{\mathbf{i}}$ can be regarded as a generalisation
of $M_{\mathbf{i}}$ for higher ranks, and it is likely that canonical basis
elements in $B_{\mathbf{i}}$ will have a common form.

Furthermore, the cones $C_{\mathbf{i}}$ play a role with respect to various
reparametrization functions of the canonical basis.
Let $\mathbf{i},\mathbf{j}$ be reduced expressions for $w_0$.
There is a reparametrization function, $R_{\mathbf{i}}^{\mathbf{j}}=
\varphi_{\mathbf{j}}\varphi_{\mathbf{i}}^{-1}:\mathbb{N}^k\rightarrow \mathbb{N}^k$
(see~\cite[\S2]{lusztig2}); a reparametrization function
$T_{\mathbf{i}}^{\mathbf{j}}=\psi_{\mathbf{j}}\psi_{\mathbf{i}}^{-1}:
Y_{\mathbf{i}} \rightarrow Y_{\mathbf{j}}$ (see~\cite[\S\S2.6, 2.7]{bz1}
and also the remark at the end of Section 2 of~\cite{nz1}).
We also define a reparametrization function
$S_{\mathbf{i}}^{\mathbf{j}}=\varphi_{\mathbf{j}}\psi_{\mathbf{i}}^{-1}:
Y_{\mathbf{i}} \rightarrow {\mathbf{N}}^k$
linking the approaches of Lusztig and Kashiwara.
The functions $R_{\mathbf{i}}^{\mathbf{j}}$ and $T_{\mathbf{i}}^{\mathbf{j}}$
are known to be piecewise-linear, and it follows that
$S_{\mathbf{i}}^{\mathbf{j}}$
is also piecewise-linear. These functions are very difficult to understand;
they have recently been described by Zelevinsky in~\cite{bz3},
but their regions of linearity are only known in types $A_1$--$A_5$,
$B_2$, $C_3$, $D_4$ and $G_2$ (see~\cite{carter3} and references therein).

Suppose now that $\mathbf{g}$ is of type $A_n$.
In~\cite{me8}, it was shown that $C_{\mathbf{i}}$ could be
expressed as the set of nonnegative integral combinations of a certain set
of $k$ spanning vectors, $v_j(\mathbf{i})$, $j=1,\ldots ,n$, 
$v_P(\mathbf{i})$, $P$ a partial quiver of type $A_n$ associated to
$\mathbf{i}$ (as described in~\cite{me8}). A partial quiver of type $A_n$ is
a quiver of type $A_n$ with some directed edges replaced by undirected edges
in such a way that the subgraph of directed edges is connected.

If an element $\mathbf{a}=(a_1,a_2,\ldots ,a_k)\in C_{\mathbf{i}}$ is
written $\mathbf{a}=(a_{\alpha})_{\alpha\in \Phi^+}$, where $\Phi^+$
is the set of
positive roots of $\mathbf{g}$ and $a_{\alpha}=a_j$ where $\alpha$ is the
$j$th root in the ordering on $\Phi^+$ induced by $\mathbf{i}$,
then the $v_j(\mathbf{i})$ and the $v_P(\mathbf{i})$ do not depend
on $\mathbf{i}$ and are described in~\cite{me8}.
These vectors have also been studied by B\'{e}dard~\cite{bedard2}
in terms of the representation theory of quivers. Here we show that if
$\mathbf{a}=(a_1,a_2,\ldots ,a_k)$ is one of these vectors then the
corresponding monomial of Kashiwara's operators $\TF_{i_1}^{a_1}
\TF_{i_2}^{a_2}\cdots \TF_{i_k}^{a_k}\cdot 1$ can be described using a
remarkable rectangle combinatorics.

For $j=1,2,\ldots ,n$ and $P$ a partial quiver associated to $\mathbf{i}$,
let $b_j(\mathbf{i})=\psi_{\mathbf{i}}^{-1}(v_j(\mathbf{i}))$ and
$b_P(\mathbf{i})=\psi_{\mathbf{i}}^{-1}(v_P(\mathbf{i}))$ be the corresponding
canonical basis elements. We will show that these elements are independent
of the reduced expression $\mathbf{i}$. Thus the canonical basis elements
arising in this way for at least one reduced expression $\mathbf{i}$ can be
labelled $b_j$, $j=1,2,\ldots ,n$, $b_P$, $P$ a partial quiver of type $A_n$.

The following conjecture, which will be shown in type $A_4$ in~\cite{me10},
suggests that these elements are important in understanding the canonical
basis. For a reduced expression $\mathbf{j}$ for $w_0$, let
$c_j(\mathbf{j})=\phi_{\mathbf{j}}(b_j)$ and let $c_P(\mathbf{j})=
\phi_{\mathbf{j}}(b_P)$ be the parametrizations of these canonical basis
elements arising from Lusztig's approach.

\begin{conjecture} (Carter and Marsh) \\
Let $\mathbf{i}$ be a reduced expression for $w_0$ and let $\mathbf{k}=
(1,3,5,\ldots 2,4,6,\ldots )$ and\linebreak
$\mathbf{k'}=(2,4,6,\ldots ,1,3,5,\ldots )$,
the opposite reduced expressions for $w_0$ arising in~\cite{lusztig2}.
Then the vectors $c_j(\mathbf{k})$,
$j=1,2,\ldots ,n$, $c_P(\mathbf{k})$, $P$ a partial quiver associated with
$\mathbf{i}$ as above form a simplicial region of linearity $X_{\mathbf{i}}$
of Lusztig's piecewise-linear function $R_{\mathbf{k}}^{\mathbf{k'}}$.
Furthermore, all such simplicial regions arise in this way, and the map
$\mathbf{i}\mapsto X_{\mathbf{i}}$ from the graph of commutation classes of
reduced expressions for $w_0$ (with edges given by the long braid relation)
to the graph of simplicial regions of linearity of $R$ (with edges given by
adjacency) is a graph isomorphism.
\end{conjecture}

In this paper we shall use the above rectangle combinatorics to
calculate the vectors $c_j(\mathbf{j})$,
$c_P(\mathbf{j})$ for the reduced expression
$\mathbf{j}=(n,n-1,n,n-2,n-1,n,\ldots ,1,2,\ldots ,n)$
for $w_0$. This will be used in a future paper of
Carter in order to calculate the vectors $c_j(\mathbf{k})$, $c_P(\mathbf{k})$,
and thus the conjectural simplicial regions of linearity of
$R_{\mathbf{k}}^{\mathbf{k'}}$. It is also possible that the above conjecture
is true for {\em any} pair of opposite reduced expressions
$\mathbf{k},\mathbf{k'}$
for $w_0$; in particular for $\mathbf{j}$ and its opposite $\mathbf{j'}$,
i.e. that the vectors $c_j(\mathbf{j})$ and $c_P^{\mathbf{j}}$ span simplicial
regions of linearity of $R_{\mathbf{j}}^{\mathbf{j'}}$.

\section{Parametrizations of the canonical basis}

For positive integers $p<q$ we denote by $[p,q]$ the set
$\{p,p+1,\ldots ,q\}$, and for a rational number $x$ we denote by
$\lceil x \rceil$ the smallest element of $\{y\in \mathbb{Z}\,:\,x\leq y\}$.

Let $\mathbf{g}$ be the simple Lie algebra over $\C$ of type $A_n$ and $U$ be the
quantized enveloping algebra of $\mathbf{g}$. Then $U$ is a $\Q(v)$-algebra
generated by the elements $E_i$, $F_i$, $K_{\mu}$, $i\in \{1,2,\ldots ,n\}$,
$\mu\in Q$, the root
lattice of $\mathbf{g}$. Let $U^+$ be the subalgebra generated by the $E_i$
and $U^-$ the subalgebra generated by the $F_i$.

Let $W$ be the Weyl group of $\mathbf{g}$ with Coxeter generators $s_1,s_2,
\ldots ,s_n$. It has a unique element $w_0$ of maximal length $k$. We shall
identify a reduced expression $s_{i_1}s_{i_2}\cdots s_{i_k}$ for $w_0$
with the $k$-tuple $\mathbf{i}=(i_1,i_2,\ldots ,i_k)$.
We shall denote by $\chi_n$ the set of all reduced expressions
${\bf i}=(i_1,i_2,\ldots ,i_k)$ for $w_0$.
Given $\mathbf{i},\mathbf{i'}\in \chi_n$,
we say that $\mathbf{i}\simeq \mathbf{i}'$ if there is a sequence of
commutations (of the form $s_is_j=s_js_i$ with $|i-j|>1$) which, when
applied to $\mathbf{i}$, give $\mathbf{i}'$.
This is an equivalence relation on $\chi_n$,
and the equivalence classes are called commutation classes.

For each $\mathbf{i}\in \chi_n$ there are two
parametrizations of the canonical basis $\B$ for $U^-$. The first arises from
Lusztig's approach to the canonical basis~\cite[\S14.4.6]{lusztig6}, and
the second arises from Kashiwara's approach~\cite{kash2}.

\pagebreak
\noindent {\bf Lusztig's Approach} \\
There is a $\Q$-algebra automorphism of $U$ which takes each $E_i$ to
$F_i$, $F_i$ to $E_i$, $K_{\mu}$ to $K_{-\mu}$ and $v$ to $v^{-1}$.
We use this automorphism to transfer Lusztig's
definition of the canonical basis in~\cite[\S3]{lusztig2} to $U^-$.

Let $T''_{i,-1}$ be the automorphism of $U$ as in~\cite[\S37.1.3]{lusztig6}.
For ${\mathbf{c}}\in {\N}^k$ and $\mathbf{i}\in \chi_n$, let
$$F_{\mathbf{i}}^{\mathbf{c}}:=F_{i_1}^{(c_1)}T''_{i_1,-1}(F_{i_2}^{(c_2)})\cdots
T''_{i_1,-1}T''_{i_2,-1}\cdots T''_{i_{k-1},-1}(F_{i_k}^{(c_k)}),$$ and
define $B_{\mathbf{i}}=\{F_{\mathbf{i}}^{\mathbf{c}}\,:\, {\mathbf{c}}\in
{\N}^k\}$. Then $B_{\mathbf{i}}$ is
the basis of PBW-type for $U^-$ corresponding to $\mathbf{i}$.
Let\,\, $\mathbf{\bar{\ }}\,$ be the ${\Q}$-algebra automorphism of $U$
taking $E_i$ to $E_i$, $F_i$ to $F_i$, and $K_{\mu}$ to $K_{-\mu}$,
for each $i\in [1,n]$ and $\mu\in Q$, and $v$ to $v^{-1}$.
Lusztig proves the following result in~\cite[\S\S2.3, 3.2]{lusztig2}.

\begin{theorem} (Lusztig) \\
The $\mathbb{Z}[v]$-span $\CL$ of $B_{\mathbf{i}}$ is independent of $\mathbf{i}$. 
Let $\pi:{{\CL}}\rightarrow {{\CL}}/{v{\CL}}$ be the natural
projection. The image $\pi(B_{\mathbf{i}})$ is also independent
of $\mathbf{i}$; we denote it by $B$. The restriction of $\pi$ to
${{\CL}}\cap \overline{{\CL}}$ is an isomorphism of $\mathbb{Z}$-modules
$\pi_1:{{\CL}}\cap \overline{{\CL}}\rightarrow {{\CL}}/{v{\CL}}$,
and $\B=\pi_1^{-1}(B)$ is a ${\Q}(v)$-basis of $U^-$, called
the canonical basis.
\end{theorem}

Lusztig's theorem provides us with a parametrization of ${\B}$, dependent on
$\mathbf{i}$. If $b\in {\B}$, we write $\varphi_{\mathbf{i}}(b)={\mathbf{c}}$,
where ${\mathbf{c}}\in {\N}^k$ satisfies $b\equiv F_{\mathbf{i}}^{\mathbf{c}}
\mod v{{\CL}}$. Note that $\varphi_{\mathbf{i}}$ is a bijection.

Given any pair ${\mathbf{i},\mathbf{j}}\in \chi_n$,
Lusztig defines in~\cite[\S2.6]{lusztig2} a function
$R_{\mathbf{i}}^{{\mathbf{j}}}=\varphi_{{\mathbf{j}}}
\varphi_{\mathbf{i}}^{-1}\,:\,{\N}^k \rightarrow {\N}^k$, which he shows
is piecewise-linear. He further shows that its regions of linearity are
significant for the canonical basis, in the sense that
elements $b$ of the canonical basis with $\varphi_{\mathbf{i}}(b)$ in the same
region of linearity of $R_{\mathbf{i}}^{{\mathbf{j}}}$ often have
similar form.

\noindent {\bf Kashiwara's approach}

Let $\widetilde{E}_i$ and $\widetilde{F}_i$ be the Kashiwara operators on
$U^-$ as defined in~\cite[\S3.5]{kash2}. Let ${\CA}\subseteq {\Q}(v)$ be
the subring of elements regular at $v=0$, and let ${\CL}'$ be the
${\CA}$-lattice spanned by the set $S$ of arbitrary products
$\widetilde{F}_{j_1}\widetilde{F}_{j_2}\cdots \widetilde{F}_{j_m}\cdot 1$ in $U^-$.
The following results were proved by Kashiwara in~\cite{kash2}.

\begin{theorem} \label{kashiwara} (Kashiwara) \\
(i) Let $\pi':{{\CL}'}\rightarrow {{\CL}'}/{v{{\CL}'}}$ be the natural
projection, and let $B'=\pi'(S)$. Then $B'$ is a ${\Q}$-basis of
${{\CL}'}/{v{{\CL}'}}$ (the crystal basis). \\
(ii) The operators
$\widetilde{E}_i$ and $\widetilde{F}_i$ each preserve ${\CL}'$ and thus act on
${{\CL}'}/{v{{\CL}'}}$. They satisfy $\widetilde{E}_i(B')\subseteq B'\cup\{0\}$ and
$\widetilde{F}_i(B')\subseteq B'$. For $b,b'\in B'$ we have $\widetilde{F}_ib=b'$
if and only if $\widetilde{E}_ib'=b$. \\
(iii) For each $b\in B'$, there is a unique element $\tilde{b}\in {{\CL}'}
\cap \overline{{\CL}'}$ such that $\pi'(\tilde{b})=b$. The set of elements
$\{\tilde{b}\,:\, b\in B'\}$ forms a basis of $U^-$, the {\em global
crystal basis} of $U^-$.
\end{theorem}

It was shown by Lusztig~\cite[2.3]{lusztig3} that the global crystal basis
of Kashiwara coincides with the canonical basis.
There is a parametrization of ${\B}$ arising from Kashiwara's approach,
again dependent on $\mathbf{i}\in \chi_n$.
Let $\mathbf{i}=(i_1,i_2,\ldots ,i_k)$ and $b\in B$.
Let $a_1$ be maximal such that $\widetilde{E}_{i_1}^{a_1}b\not\equiv 0 \mod
v{{\CL}'}$;
let $a_2$ be maximal such that
$\widetilde{E}_{i_2}^{a_2}\widetilde{E}_{i_1}^{a_1}b\not\equiv 0 \mod v{{\CL}'}$,
and so on, so that
$a_k$ is maximal such that $\widetilde{E}_{i_k}^{a_k}\widetilde{E}_{i_{k-1}}^{a_{k-1}}
\cdots \widetilde{E}_{i_2}^{a_2}\widetilde{E}_{i_1}^{a_1}b\not\equiv 0 \mod
v{{\CL}'}$.
We write $\psi_{\mathbf{i}}(b)=\mathbf{a}=(a_1,a_2,\ldots ,a_k)$ and
we have $b\equiv \TF_{i_1}^{a_1}\TF_{i_2}^{a_2}\cdots \TF_{i_k}^{a_k}\cdot 1
\mod v{\CL}'$.
This is the crystal string of $b$ --- see~\cite[\S2]{bz1},~\cite[\S2]{nz1}
and~\cite{kash4}.
The map $\psi_{\mathbf{i}}$ is injective (see~\cite[\S2.5]{nz1}).
Its image is a cone, which we shall call the
{\em string cone} associated to $\bf i$, $Y_{\mathbf{i}}=\psi_{\mathbf{i}}(\B)$.
We define a function which compares Kashiwara's approach with
Lusztig's approach. For $\mathbf{i}$, $\mathbf{j}\in \chi_n$,
let $S_{\mathbf{i}}^{\mathbf{j}}=\varphi_{\mathbf{j}}
\psi_{\mathbf{i}}^{-1}\,:\,Y_{\mathbf{i}}\rightarrow {\N}^k$.

\section{The Lusztig cones and their spanning vectors}

Lusztig~\cite{lusztig7} introduced certain regions which, in low rank,
give rise to canonical basis elements of a particularly simple form. 
The {\em Lusztig cone}, $C_{\mathbf{i}}$, corresponding
to $\mathbf{i}\in \chi_n$, is defined to be the set of
points ${{\mathbf{a}}}\in {\N}^k$
satisfying the following inequalities: \\
(*) For every pair $s,s'\in [1,k]$ with $s<s'$, $i_s=i_{s'}=i$ and
$i_p\not=i$ whenever $s<p<s'$, we have
$$(\sum_p a_p)-a_s - a_{s'}\geq 0,$$
where the sum is over all $p$ with $s<p<s'$ such that $i_p$ is joined to $i$
in the Dynkin diagram (we call such a pair $s,s'$ a {\em minimal pair}).


It was shown by Lusztig~\cite{lusztig7} that if $\mathbf{a}\in C_{\mathbf{i}}$
then the monomial $F_{i_1}^{(a_1)}F_{i_2}^{(a_2)}\cdots
F_{i_k}^{(a_k)}$ lies in the canonical basis ${\B}$, provided $n=1,2,3$.
The author~\cite{me7} showed that this remains true if $n=4$.
Counter-examples of Xi~\cite{xi4} and Reineke~\cite{reineke2}
show that this is not true in general.

{\bf Remark:} It can be shown that, for any $\mathbf{i}\in \chi_n$,
we have $C_{\mathbf{i}}\subseteq Y_{\mathbf{i}}$, i.e., the Lusztig cone is
contained in the Kashiwara cone. The methods of this paper could be used to
do this in type $A_n$, but the author has recently learnt that
A. Premat~\cite{premat1} has a more succinct proof for all simply-laced cases.

The reduced expression $\mathbf{i}$ defines an ordering on the set $\Phi^+$ of
positive roots of the root system associated to $W$.
We write $\alpha^j=s_{i_1}s_{i_2}\cdots s_{i_{j-1}}(\alpha_{i_j})$ for
$j=1,2,\ldots ,k$. Then $\Phi^+=\{\alpha^1,\alpha^2,\ldots \alpha^k\}$.
For ${{\mathbf{a}}}=(a_1,a_2,\ldots
,a_k)\in {\mathbb{Z}}^k$, write $a_{\alpha^j}=a_j$.
If $\alpha=\alpha_{ij}:=\alpha_i+\alpha_{i+1}+\cdots +\alpha_{j-1}$
with $i<j$, we also write $a_{ij}$ for $a_{\alpha_{ij}}$.

\begin{figure}
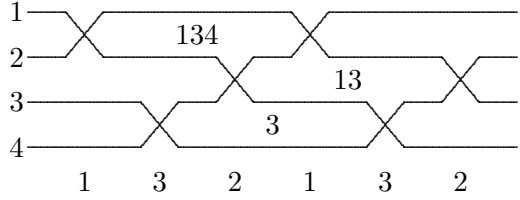

\beginpicture

\setcoordinatesystem units <0.5cm,0.3cm>             
\setplotarea x from -10 to 14, y from -3 to 8       
\linethickness=0.5pt           

\put{$1$}[c] at -0.3 6 %
\put{$2$}[c] at -0.3 4 %
\put{$3$}[c] at -0.3 2 %
\put{$4$}[c] at -0.3 0 %

\put{$1$}[c] at 1.5 -1.5 %
\put{$3$}[c] at 3.5 -1.5 %
\put{$2$}[c] at 5.5 -1.5 %
\put{$1$}[c] at 7.5 -1.5 %
\put{$3$}[c] at 9.5 -1.5 %
\put{$2$}[c] at 11.5 -1.5 %

\setlinear \plot 0 0  3 0 / %
\setlinear \plot 3 0  4 2 / %
\setlinear \plot 4 2  5 2 / %
\setlinear \plot 5 2  6 4 / %
\setlinear \plot 6 4  7 4 / %
\setlinear \plot 7 4  8 6 / %
\setlinear \plot 8 6  13 6 / %

\setlinear \plot 0 2  3 2 / %
\setlinear \plot 3 2  4 0 / %
\setlinear \plot 4 0  9 0 / %
\setlinear \plot 9 0  10 2 / %
\setlinear \plot 10 2  11 2 / %
\setlinear \plot 11 2  12 4 / %
\setlinear \plot 12 4  13 4 / %

\setlinear \plot 0 4  1 4 / %
\setlinear \plot 1 4  2 6 / %
\setlinear \plot 2 6  7 6 / %
\setlinear \plot 7 6  8 4 / %
\setlinear \plot 8 4  11 4 / %
\setlinear \plot 11 4  12 2 / %
\setlinear \plot 12 2  13 2 / %

\setlinear \plot 0 6  1 6 / %
\setlinear \plot 1 6  2 4 / %
\setlinear \plot 2 4  5 4 / %
\setlinear \plot 5 4  6 2 / %
\setlinear \plot 6 2  9 2 / %
\setlinear \plot 9 2  10 0 / %
\setlinear \plot 10 0  13 0 / %

\put{$134$}[c] at 4.5 5 %
\put{$13$}[c] at 8.5 3 %
\put{$3$}[c] at 6.5 1 %
\endpicture
\caption{The Chamber Diagram of $(1,3,2,1,3,2)$.\label{chamberdiagram}}
\end{figure}

We shall need the chamber diagram (chamber ansatz) for $\mathbf{i}$
defined in~\cite[\S\S1.4, 2.3]{bfz1}.
We take $n+1$ strings, numbered from top to bottom, and write 
$\mathbf{i}$ from left to right along the bottom of the diagram.
Above a letter $i_j$ in $\mathbf{i}$,
the $i_j$th and $(i_j+1)$st strings from the top above $i_j$ cross.
Thus, for example, in the case $n=3$ with $\mathbf{i}=(1,3,2,1,3,2)$, the
chamber diagram is shown in Figure~\ref{chamberdiagram}.

We denote the chamber diagram of $\mathbf{i}$ by CD($\mathbf{i}$).
A {\em chamber} will be defined as a pair $(c,{\mathbf{i}})$, where $c$ is a
bounded component of the complement of CD($\mathbf{i}$).
Each chamber $(c,{\mathbf{i}})$ can be labelled with the numbers of the
strings passing below it, denoted $l(c,{\mathbf{i}})$.
Following~\cite{bfz1}, we call such a label a chamber set.
For example, the chamber sets corresponding to the $3$ bounded chambers in
Figure~\ref{chamberdiagram} are $134,3,13$. Note that the set of chamber
sets of $\bf i$ is independent of its commutation class
(so we can talk also of the chamber sets of such a commutation class).
Every subset of $[1,n+1]$ can arise as a chamber set for some $\bf i$,
except for subsets of the form $[1,j]$ and $[j,n+1]$ for $1\leq j\leq n+1$;
this is observed in the proof of~\cite[Theorem 2.7.1]{bfz1}.

We now recall some of the results from~\cite{me8} that we shall need.
We have by~\cite[\S4]{me8} that there is a matrix $P_{\mathbf{i}}\in
GL_k(\mathbb{Z})$ such that
$C_{\bf i}=\{{\bf a}\in \mathbb{Z}^k\,:\,P_{\bf i}{\bf a}\geq 0\},$
where, for ${\bf z}\in \mathbb{Z}^k$, $\mathbf{z}\geq 0$ means that each entry in
$\bf z$ is nonnegative.
It follows that $C_{\bf i}$ is the set of nonnegative linear combinations of
the columns of $Q_{\bf i}=P_{\mathbf{i}}^{-1}$. These columns were described
in~\cite{me8}, using the concept of a {\em partial quiver}.

A partial quiver $P$ of type $A_n$ is
a quiver of type $A_n$ which has some (or none) of its arrows replaced
by undirected edges, in such a way that the subgraph obtained by deleting
undirected edges and any vertex incident only with undirected edges is
connected.
We number the edges of a partial quiver from $2$ to $n$, starting at
the right hand end.
We write $P$ as a sequence of $n-1$ symbols, $L$, $R$ or $-$,
where $L$ denotes a leftward arrow, $R$ a rightward arrow, and $-$ an
undirected edge. Thus any partial quiver will be of form $---***---$, where
the $*$'s denote $L$'s or $R$'s.

If $P,P'$ are partial quivers we write $P'\geq P$ (or $P\leq P'$)
if every edge which is directed in $P$ is directed in $P'$ and is oriented in the
same way, and say that $P$ is a sub partial quiver of $P'$.
For example, $---LRLL-\ \leq RLRLRLLL$.

It is known that if $\bf i$ is compatible with a quiver $Q$
(in the sense of~\cite[\S4.7]{lusztig2}) then the
set of reduced expressions for $w_0$ compatible with $Q$ is precisely
the commutation class of $\chi_n$ containing $\bf i$.
We say that this commutation class is compatible with $Q$.
Berenstein, Fomin and Zelevinsky (see~\cite[\S4.4]{bfz1})
describe a method for constructing $\bf i\in \chi_n$ compatible with any given
quiver $Q$, as follows.

\begin{figure}
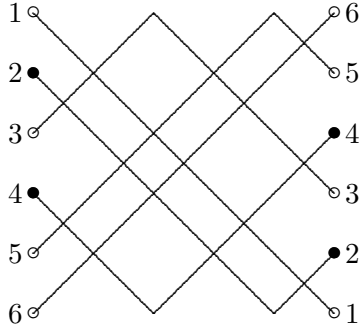

\begin{center}
\beginpicture
\setcoordinatesystem units <0.8cm,0.8cm>             
\setplotarea x from -7 to 7, y from 0 to 7       

\linethickness=0.5pt           

\multiput{$\circ$} at 1 1 *5 0 1 /
\multiput{$\circ$} at 6 1 *5 0 1 /

\put{$\bullet$} at 1 3 %
\put{$\bullet$} at 1 5 %
\put{$\bullet$} at 6 2 %
\put{$\bullet$} at 6 4 %

\put{$1$}[c] at 0.7 6 %
\put{$2$}[c] at 0.7 5 %
\put{$3$}[c] at 0.7 4 %
\put{$4$}[c] at 0.7 3 %
\put{$5$}[c] at 0.7 2 %
\put{$6$}[c] at 0.7 1 %

\put{$6$}[c] at 6.3 6 %
\put{$5$}[c] at 6.3 5 %
\put{$4$}[c] at 6.3 4 %
\put{$3$}[c] at 6.3 3 %
\put{$2$}[c] at 6.3 2 %
\put{$1$}[c] at 6.3 1 %

\setlinear \plot 1 6    6 1 / %

\setlinear \plot 1 5    5 1 / %
\setlinear \plot 5 1    6 2 / %

\setlinear \plot 1 4    3 6 / %
\setlinear \plot 3 6    6 3 / %

\setlinear \plot 1 3    3 1 / %
\setlinear \plot 3 1    6 4 / %

\setlinear \plot 1 2    5 6 / %
\setlinear \plot 5 6    6 5 / %

\setlinear \plot 1 1    6 6 / %

\endpicture
\caption{The arrangement Arr$\{2,4\}$.\label{arrangement}}
\end{center}
\end{figure}

Suppose $Q$ is a quiver of type $A_n$. Let $\Lambda\subseteq [2,n]$ be the
set of all edges of $Q$ pointing to the left. Berenstein, Fomin and
Zelevinsky construct an arrangement Arr$(\Lambda)$.
Consider a square in the plane, with horizontal
and vertical sides.
Put $n+1$ points onto the left-hand edge of the square, equally spaced,
numbered $1$ to $n+1$ from top to bottom, and
do the same for the right-hand edge, but number the points
from bottom to top. ${\rm Line}_h$ joins point $h$ on the left with
point $h$ on the right. For $h=1,n+1$, ${\rm Line}_h$ is a diagonal
of the square. For $h\in [2,n]$, ${\rm Line}_h$ is a union of two
line segments of slopes $\pi /4$ and $-\pi /4$. There are precisely
two possibilities for $\mbox{Line}_h$. If $h\in \Lambda$, the
left segment has slope $-\pi/4$, while the right one has slope $\pi/4$;
for $h\in [2,n]\setminus \Lambda$, it goes the other way round.
Berenstein, Fomin and Zelevinsky give the example of the arrangement 
for $n=5$ and $\Lambda=\{2,4\}$, which we show in Figure~\ref{arrangement} (points
corresponding to elements of $\Lambda$ are indicated by filled circles).
They show that $\bf i\in \chi_n$ is compatible with the quiver $Q$
if and only if the chamber diagram CD($\bf i$) is isotopic to Arr$(\Lambda)$.

If $P$ is a partial quiver, let $l(P)$ be the subset of $[1,n+1]$
defined in the following way.
Put $l_1(P)=\{j\in [2,n]\,:\,\mbox{edge\ $j$\ of\ $P$\ is\ an\ 
$L$}\}$.
If the rightmost directed edge of $P$ is an $R$, and this is
in position $i$, then let $l_2(P)=[1,i-1]$, otherwise the empty set.
If the leftmost directed edge of $P$ is an $R$, and this is
in position $j$, then let $l_3(P)=[j+1,n+1]$, otherwise the empty set.
Then put $l(P)=l_1(P)\cup l_2(P)\cup l_3(P)$.
We know from~\cite[\S5.4]{me8}
that $l$ sets up a one-to-one correspondence between
chamber sets and partial quivers; for a chamber $(c,{\bf i})$, we denote
by $P(c,{\bf i})$ the partial quiver corresponding to $l(c,{\bf i})$.
Also, if $\bf i$ is compatible with a quiver $Q$, then the partial quivers
corresponding to its chambers are precisely those $P\leq Q$.

Each minimal pair in ${\bf i}$ corresponds naturally to a chamber
$(c,{\bf i})$. Thus, given a chamber $(c,{\bf i})$,
there is a corresponding row of $P_{\bf i}$ and
therefore a corresponding column of $Q_{\bf i}$, that is,
a spanning vector of $C_{\bf i}$. We denote this spanning vector by
$a(c,{\bf i})$.
Similarly, the other $n$ rows of $P_{\bf i}$ correspond to inequalities
of the form $a_{\alpha_j}\geq 0$ for $\alpha_j$ a simple root, $1\leq j\leq
n$. We denote the corresponding spanning vectors by $a(j,{\bf i})$.

We consider $a(c,{\bf i})=(a_1,a_2,\ldots ,a_k)$ as a multiset $M(c,{\bf i})$ of
positive roots, where each positive root $\alpha_{ij}$
occurs with multiplicity, $a_{ij}$.
Similarly we write $M(j,{\bf i})$ for the multiset corresponding to 
$a(j,{\bf i})$, for $j=1,2,\ldots ,k$. These multisets have been described
in~\cite{me8}.

\begin{theorem} (See~\cite[\S5.12]{me8}). \label{multiset} \ \\
(a) The multiset spanning vector $M(c,{\bf i})$ depends only upon
$P(c,{\bf i})$. For a partial quiver $P$ we choose a chamber
$(c,{\bf i})$ such that $P(c,{\bf i})=P$ and write $M(P)=M(c,{\bf i})$. \\
(b) The multiset spanning vector $M(j,{\bf i})$ depends only upon $j$;
we denote it by $M(j)$. \\
(c) For a partial quiver $P$, $M(P)$ is given by the following construction.
We say that a sub partial quiver $Y$ of $P$ is a {\em component}
of $P$ if all of its edges are oriented in the same way and it is
maximal in length with this property.
We say $Y$ has type $L$ if its edges are oriented to the left, and type $R$
if its edges are oriented to the right.
For such a component $Y$ of $P$, let
$a(Y)$ be the number of the leftmost
edge to the right of the component, and let $b(Y)$ be the number
of the rightmost edge to the left of the component.
Let $M(Y)$ be the set
of positive roots $\alpha_{ij}$ satisfying $1\leq i\leq a(Y)\leq b(Y)\leq j
\leq n+1$. Then $M(P)$ is the multiset union (i.e.\ adding
multiplicities) of the $M(Y)$, where $Y$ varies over all components
of $P$, with the multiplicity $m_{ij}$ of $\alpha_{ij}$ replaced by
$n_{ij}=\lceil\frac{1}{2}m_{ij}\rceil$.
Furthermore, for $j=1,2,\ldots n$, $M(j)$ is the set of positive roots
$\alpha_{pq}$ satisfying $1\leq p\leq j\leq j+1\leq q\leq n+1$.
\end{theorem}

\section{Canonical basis elements corresponding to spanning vectors}

Recall that for $\mathbf{i},\mathbf{i'}\in \chi_n$, the map
$T_{\bf i}^{\bf i'}$ is defined by
$T_{\bf i}^{\bf i'}({\bf a})=\psi_{\bf i'}\psi_{\bf i}^{-1}({\bf a})$, for
${\bf a}\in Y_{\bf i}$.
Define $T_2:\,\mathbb{Z}^2\rightarrow \mathbb{Z}^2$ by $T_2(a,b):=(b,a)$, and
$T_3:\,\mathbb{Z}^3\rightarrow \mathbb{Z}^3$ by $T_3(a,b,c):=(\max(c,b-a),a+c,
\min(a,b-c))$. The maps $T_{\bf i}^{\bf i'}$ satisfy the following:

\begin{prop} \label{Tdescription} (Berenstein and Zelevinsky) \\
(a) Suppose ${\bf i,i'}\in \chi_n$ and that $\mathbf{i}=\mathbf{i'}$ except
that $(i_s,i_{s+1})=(i,j)$ and $(i'_s,i'_{s+1})=(j,i)$ for some index $s$ and
$i,j\in [1,n]$ with $|i-j|>1$. Let $\mathbf{a}\in \mathbb{N}^k$ and
write ${\bf a'}=T_{\bf i}^{\bf i'}({\bf a})$. Then
we have $a'_r=a_r$ for $r<s$ or $r>s+1$ and
$(a'_s,a'_{s+1})=T_2(a_s,a_{s+1})$. \\
(b) Suppose ${\bf i,i'}\in \chi_n$, and that $\mathbf{i}=\mathbf{i'}$ except
that $(i_s,i_{s+1},i_{s+2})=(i,j,i)$ and $(i'_s,i'_{s+1},i'_{s+2})=(j,i,j)$ for
some index $s$ and $i,j\in [1,n]$ with $|i-j|=1$.
Let $\mathbf{a}\in \mathbb{N}^k$ and write
${\bf a'}=T_{\bf i}^{\bf i'}({\bf a})$. Then we have $a'_r=a_r$ for $r<s$
or $r>s+2$ and $(a'_s,a'_{s+1},a'_{s+2})=T_3(a_s,a_{s+1},a_{s+2})$.
\end{prop}

\noindent {\bf Proof}: See~\cite[\S\S2.6, 2.7]{bz1}.~$\square$

We need the following Lemma, which will help us to understand how spanning
vectors of Lusztig cones are related.

\begin{lemma} \label{shortT}
(a) Suppose ${\bf i,i'}\in \chi_n$ and that $\mathbf{i}=\mathbf{i'}$
except that $(i_s,i_{s+1})=(i,j)$ and $(i'_s,i'_{s+1})=(j,i)$ for some index
$s$ and $i,j\in [1,n]$ with $|i-j|>1$.
Let $\mathbf{a}=a(c,\mathbf{i})$, respectively $a(j,\mathbf{i})$, be a
spanning vector of
$C_{\mathbf{i}}$, and let $\mathbf{a'}=a(c',\mathbf{i'})$, respectively
$\mathbf(j,\mathbf{i'})$,
be the corresponding spanning vector of $C_{\mathbf{i'}}$.
Here $(c',\mathbf{i'})$ is the chamber that $(c,\mathbf{i})$ becomes when the
commutation relation
$(i,j)\rightarrow (j,i)$ is applied to $\mathbf{i}$. Then $a'_r=a_r$ for all
$r<s$ or $r>s+1$, while $(a'_s,a'_{s+1})=T_2(a_s,a_{s+1})$.\\
(b) Suppose ${\bf i,i'}\in \chi_n$ and that $\mathbf{i}=\mathbf{i'}$
except that $(i_s,i_{s+1},i_{s+2})=(i,j,i)$ and $(i'_s,i'_{s+1},i'_{s+2})=
(j,i,j)$ for some index $s$ and $i,j\in [1,n]$ with $|i-j|=1$.
Let $\mathbf{a}=a(c,\mathbf{i})$, respectively $a(j,\mathbf{i})$, be a
spanning vector of $C_{\mathbf{i}}$.
In this case, we suppose further that $(c,\mathbf{i})$ is not the chamber of
$CD(\mathbf{i})$ involving the substring $(i,j,i)$, so that there is a chamber
$(c',\mathbf{i'})$ with $P(c,\mathbf{i})=P(c',\mathbf{i'})$.
Let $\mathbf{a'}=a(c',\mathbf{i'})$, respectively $\mathbf(j,\mathbf{i'})$.
Then $a'_r=a_r$ for all $r<s$ or $r>s+1$, while $(a'_s,a'_{s+1},a'_{s+2})=
T_3(a_s,a_{s+1},a_{s+2})$.
\end{lemma}

\noindent {\bf Proof}: Recall that the spanning vectors of $C_{\bf i}$ are
obtained by inverting the defining matrix $P_{\bf i}$ of $C_{\bf i}$. 
Part (a) follows immediately from the definition of $P_{\bf i}$.
Suppose we are in situation (b). Consider the case where
${\bf a}=a(c,{\bf i})$. Since $(c,{\bf i})$ is not the
chamber corresponding to $(i,j,i)=(i_s,i_{s+1},i_{s+2})$, we must have that
${\bf a}$ is perpendicular to the row of $P_{\bf i}$ corresponding to this
chamber. In the case where ${\bf a}=a(j,{\bf i})$, this is also true.
So, in either case, we see that $a_{s+1}=a_s+a_{s+2}$. It follows from the
definition of $T_3$ that
$T_3(a_s,a_{s+1},a_{s+2})=(a_{s+2},a_{s+1},a_s)$. Since the ordering on the
positive roots induced by $\bf i$ agrees with that induced by $\bf i'$
except that the roots in positions $s$ and $s+2$ are exchanged, the fact
that $\bf a'$ is as claimed in the statement above follows from
Theorem~\ref{multiset} (a) and (b).~$\square$

Using Proposition~\ref{Tdescription} and Lemma~\ref{shortT} we have:

\begin{cor}
Let ${\bf i,i'}\in \chi_n$, and let $(c,{\bf i}),(c',{\bf i'})$ be chambers
with $P(c,{\bf i})=P(c',{\bf i'})$. Then $T_{\bf i}^{\bf i'}(a(c,{\bf i}))=
a(c',{\bf i'})$. Also, for $j\in [1,n]$, $T_{\bf i}^{\bf i'}(a(j,{\bf i}))=
a(j,{\bf i'})$.
\end{cor}

We therefore have:

\begin{theorem} \label{natural}
Let $\mathbf{i}\in \chi_n$.
For $j=1,2,\ldots ,n$ and $P$ a partial quiver associated to $\mathbf{i}$,
let $b_j(\mathbf{i})=\psi_{\mathbf{i}}^{-1}(v_j(\mathbf{i}))$ and
$b_P(\mathbf{i})=\psi_{\mathbf{i}}^{-1}(v_P(\mathbf{i}))$ be the corresponding
canonical basis elements. Then these elements are independent
of the reduced expression $\mathbf{i}$. Thus the canonical basis elements
arising in this way for at least one reduced expression $\mathbf{i}$ can be
labelled $b_j$, $j=1,2,\ldots ,n$, $b_P$, $P$ a partial quiver of type $A_n$.
\end{theorem}

\section{Rectangles} \label{rectangles}
The above provides us with a description of the spanning vectors of a
Lusztig cone as multisets of positive roots. We would also like a description
of the monomial of Kashiwara operators corresponding to such a vector.
We will use this in the next section to calculate Lusztig parametrizations
$c_j(\mathbf{j})$, $c_P(\mathbf{j})$ of the canonical basis elements
arising in Theorem~\ref{natural}.
We shall show that these monomials can be described by a remarkable `rectangle
combinatorics' (invented by R. W. Carter) which we shall now describe.

\begin{definition} \rm \label{diagonal}
Let $P$ be a partial quiver, and $Y$ a component of $P$.
Let $a=a(Y)$ and let $b=b(Y)$ (see Theorem~\ref{multiset}).
We define the {\em rectangle} $\rho(Y)$ of $Y$, as follows. If $Y$ is of type
$L$, then $\rho(Y)$ is as in the left hand diagram in Figure~\ref{rectanglesf};
if $Y$ is of type $R$ then $\rho(Y)$ is as in the right hand diagram in
Figure~\ref{rectanglesf}. We also define, for $j=1,2,\ldots ,n$,
$\rho(j)$ to be the type $L$ rectangle with $a=j$
and $b=j+1$ (for these values of $a$ and $b$, both diagrams in
Figure~\ref{rectanglesf} coincide).

\begin{figure}
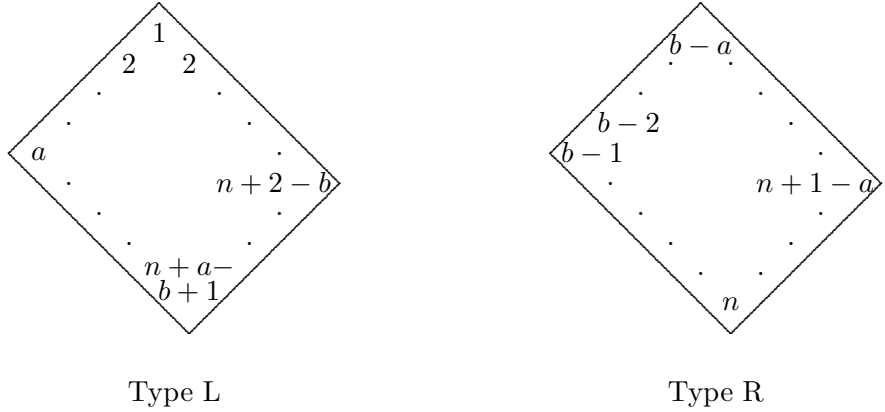

\beginpicture

\setcoordinatesystem units <0.8cm,0.8cm>             
\setplotarea x from -4 to 14, y from -4.5 to 3       

\linethickness=0.5pt           


\put{$a$}[c] at 0.5 0 %
\put{$.$}[c] at 1 0.5 %
\put{$.$}[c] at 1.5 1.0 %
\put{$2$}[c] at 2 1.5 %
\put{$1$}[c] at 2.5 2.0 %

\put{$2$}[c] at 3 1.5 %
\put{$.$}[c] at 3.5 1 %
\put{$.$}[c] at 4 0.5 %
\put{$.$}[c] at 4.5 0 %
\put{$n+2-b$}[c] at 4.4 -0.5 %

\put{$.$}[c] at 4.5 -1 %
\put{$.$}[c] at 4 -1.5 %
\put{$n+a-$}[c] at 3 -1.9 %
\put{$b+1$}[c] at 3 -2.3 %

\put{$.$}[c] at 2 -1.5 %
\put{$.$}[c] at 1.5 -1 %
\put{$.$}[c] at 1 -0.5 %

\setlinear \plot 0 0  2.5 2.5 / %
\setlinear \plot 2.5 2.5  5.5 -0.5 / %
\setlinear \plot 5.5 -0.5 3 -3 / %
\setlinear \plot 3 -3  0 0 / %

\put{Type L}[c] at 2.75 -4 %

\put{$b-1$}[c] at 9.7 0 %
\put{$b-2$}[c] at 10.3 0.5 %
\put{$.$}[c] at 10.5 1.0 %
\put{$.$}[c] at 11 1.5 %
\put{$b-a$}[c] at 11.5 1.8 %

\put{$.$}[c] at 12 1.5 %
\put{$.$}[c] at 12.5 1 %
\put{$.$}[c] at 13 0.5 %
\put{$.$}[c] at 13.5 0 %
\put{$n+1-a$}[c] at 13.4 -0.5 %

\put{$.$}[c] at 13.5 -1 %
\put{$.$}[c] at 13 -1.5 %
\put{$.$}[c] at 12.5 -2 %
\put{$n$}[c] at 12 -2.5 %

\put{$.$}[c] at 11.5 -2 %
\put{$.$}[c] at 11 -1.5 %
\put{$.$}[c] at 10.5 -1 %
\put{$.$}[c] at 10 -0.5 %

\setlinear \plot 9 0  11.5 2.5 / %
\setlinear \plot 11.5 2.5  14.5 -0.5 / %
\setlinear \plot 14.5 -0.5 12 -3 / %
\setlinear \plot 12 -3  9 0 / %

\put{Type R}[c] at 11.75 -4 %

\endpicture

\caption{The rectangles $\rho(Y)$.\label{rectanglesf}}
\end{figure}

In each case, the rectangle is filled in with numbers in diagonal rows,
with the numbers along a diagonal row from top right to bottom left
increasing by 1 at each step and similarly from top left to bottom right.
We define the diagram $D(P)$ of $P$ in the following way. Go through the
components $Y$ of $P$ one by one, from left to right.
Fit the corresponding
rectangles $\rho(Y)$ together, so that if a component of
type $L$
is followed by one of type $R$, the rectangles share leftmost corners, and
if one of type $R$ is followed by one of type $L$, they share rightmost
corners. In each case, it is easy to see that the overlapping numbers
agree. We associate a multiplicity $m(s)$ to each position in $D(P)$ where
a number $s$ is placed; if $t$ is the number of rectangles $\rho(Y)$ that $s$
lies in, then $m(s)$ is defined to be $\lceil t/2\rceil$.
We label each rectangle with the corresponding component
$Y$ of $P$; on the left hand corner if $Y$ is of type $R$, and on the right
hand corner if $Y$ is of type $L$. We define, for $j=1,2,\ldots ,n$, $D(j)$
to be the diagram $\rho(j)$, with each position given multiplicity $1$.
\end{definition}

\begin{figure}
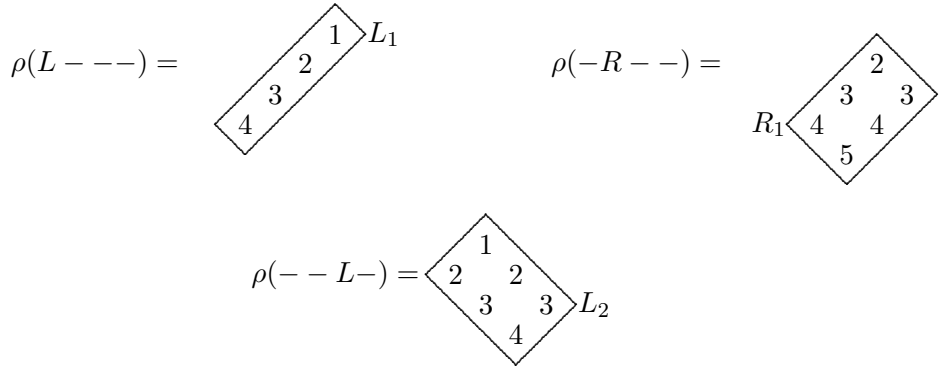

\beginpicture

\setcoordinatesystem units <0.8cm,0.8cm>             
\setplotarea x from -6 to 10, y from -4.5 to 3       

\linethickness=0.5pt           


\put{$\rho(L---)=$} at -2 1 %

\put{$4$}[c] at 0.5 0 %
\put{$3$}[c] at 1 0.5 %
\put{$2$}[c] at 1.5 1.0 %
\put{$1$}[c] at 2 1.5 %

\setlinear \plot 0 0  2 2 / %
\setlinear \plot 2 2  2.5 1.5 / %
\setlinear \plot 2.5 1.5  0.5 -0.5 / %
\setlinear \plot 0.5 -0.5  0 0 / %

\put{$L_1$} at 2.8 1.5 %


\put{$\rho(-R--)=$} at 7 1 %

\put{$4$}[c] at 10 0 %
\put{$3$}[c] at 10.5 0.5 %
\put{$2$}[c] at 11 1.0 %

\put{$5$}[c] at 10.5 -0.5 %
\put{$4$}[c] at 11 0 %
\put{$3$}[c] at 11.5 0.5 %

\setlinear \plot 9.5 0  11 1.5 / %
\setlinear \plot 11 1.5  12 0.5 / %
\setlinear \plot 12 0.5  10.5 -1 / %
\setlinear \plot 10.5 -1  9.5 0 / %

\put{$R_1$} at 9.2 0 %


\put{$\rho(--L-)=$} at 2 -2.5 %

\put{$2$}[c] at 4 -2.5 %
\put{$1$}[c] at 4.5 -2 %

\put{$3$}[c] at 4.5 -3 %
\put{$2$}[c] at 5 -2.5 %

\put{$4$}[c] at 5 -3.5 %
\put{$3$}[c] at 5.5 -3 %

\setlinear \plot 3.5 -2.5  4.5 -1.5 / %
\setlinear \plot 4.5 -1.5  6 -3 / %
\setlinear \plot 6 -3  5 -4 / %
\setlinear \plot 5 -4  3.5 -2.5 / %

\put{$L_2$} at 6.3 -3 %

\endpicture
\caption{Examples of the rectangles $\rho(Y)$.\label{egrectangles}}
\end{figure}

\noindent {\bf Example} \\
We consider the case when $n=5$ and $P=LRL-$. Then $P$ has three components,
$L_1=L---$, $R_1=-R--$ and $L_2=--L-$.
The diagrams $\rho(Y)$ are as in Figure~\ref{egrectangles}.

\begin{figure}
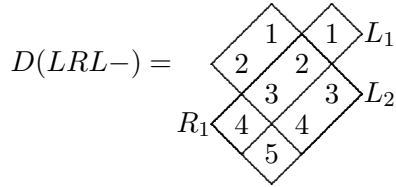

\beginpicture

\setcoordinatesystem units <0.8cm,0.8cm>             
\setplotarea x from -9 to 3, y from -1 to 2.75       

\linethickness=0.5pt           

\put{$D(LRL-)=$} at -2 1 %


\put{$4$}[c] at 0.5 0 %
\put{$3$}[c] at 1 0.5 %
\put{$2$}[c] at 1.5 1.0 %
\put{$1$}[c] at 2 1.5 %

\setlinear \plot 0 0  2 2 / %
\setlinear \plot 2 2  2.5 1.5 / %
\setlinear \plot 2.5 1.5  0.5 -0.5 / %
\setlinear \plot 0.5 -0.5  0 0 / %


\put{$5$}[c] at 1 -0.5 %
\put{$4$}[c] at 1.5 0 %
\put{$3$}[c] at 2 0.5 %

\setlinear \plot 0 0  1.5 1.5 / %
\setlinear \plot 1.5 1.5  2.5 0.5 / %
\setlinear \plot 2.5 0.5  1 -1 / %
\setlinear \plot 1 -1  0 0 / %


\put{$2$}[c] at 0.5 1 %
\put{$1$}[c] at 1 1.5 %

\setlinear \plot 0 1  1 2 / %
\setlinear \plot 1 2  2.5 0.5 / %
\setlinear \plot 2.5 0.5  1.5 -0.5 / %
\setlinear \plot 1.5 -0.5  0 1 / %

\put{$L_1$}[c] at 2.8 1.5 %
\put{$R_1$}[c] at -0.3 0 %
\put{$L_2$}[c] at 2.8 0.5 %


\endpicture
\caption{The diagram $D(LRL-)$.\label{lrl}}
\end{figure}

We put them together, so that the $4$ in $\rho(L---)$ matches the leftmost
$4$ in $\rho(-R--)$ and so that the rightmost $3$ in $\rho(-R--)$ matches
the rightmost $3$ in $\rho(--L-)$, to get the complete diagram $D(LRL-)$;
see Figure~\ref{lrl}. The numbers $2$ and $3$ in the middle of the diagram
each have multiplicity $2$.

\begin{definition} \rm \label{diagonalrow}
Suppose $P$ is a partial quiver or $P=j$ for some $j\in [1,n]$.
We define a {\em diagonal row} of $D(P)$ to be a diagonal row of numbers in $D(P)$
running from top left to bottom right.
Such a diagram $D(P)$ defines a sequence $\mu(P)$ which is obtained by
reading off the digits in the diagonal rows in the diagram, starting with
the bottom left row. Each digit is repeated according to its multiplicity.
In our example, $\mu(LRL-)=(4,5,2,3,3,4,1,2,2,3,1)$.
There is a corresponding monomial $F(P)$ in $U^-$, obtained by
replacing each digit, $s$, occurring with multiplicity $m(s)$, by the
divided power $F_s^{(m(s))}$. In our example, $F(LRL-)=
F_4F_5F_2F_3^{(2)}F_4F_1F_2^{(2)}F_3F_1$.
Similarly, there is a corresponding monomial in the Kashiwara operators
which we denote $\widetilde{F}(P)$. So in our example,
$\widetilde{F}(LRL-)=\widetilde{F}_4\widetilde{F}_5\widetilde{F}_2\widetilde{F}_3^{2}
\widetilde{F}_4\widetilde{F}_1\widetilde{F}_2^{2}\widetilde{F}_3\widetilde{F}_1\cdot 1$.
\end{definition}
The following theorem was conjectured by Carter:

\begin{theorem} \label{spanningmonomial}
Suppose $\bf i\in \chi_n$ is compatible with a quiver $Q$,
and let $(c,{\bf i})$ be a chamber in CD$({\bf i})$.
Let $P=P(c,{\bf i})$ be the corresponding partial quiver (so $P\leq Q$),
with ${\bf a}=a(c,{\bf i})=(a_1,a_2,\ldots ,a_k)$.
Alternatively, let $P=j\in [1,n]$. In either case,
$$F_{i_1}^{(a_1)}F_{i_2}^{(a_2)}\cdots F_{i_k}^{(a_k)}=F(P).$$
In particular, the left hand side depends only upon $P$.
\end{theorem}

\noindent {\bf Proof:}
Suppose first that $P$ is a partial quiver, and let $Y$ be a component of $P$.
Our first step will be to understand where the roots in $M(Y)$ occur in
CD$({\bf i})$ (as labels of crossings).
Let us suppose in the first case that $Y$ is of type $L$,
with $a=a(Y)$ and $b=b(Y)$.
We draw the diagram for the arrangement corresponding to $Q$, defined
by Berenstein, Fomin and Zelevinsky.
We mark the crossing of strings $i$ and $j$, $i\leq j$
by the root $\alpha_{ij}=\alpha_i+\alpha_{i+1}+\cdots +\alpha_{j-1}$, and
we write the quiver $Q$ up the left hand side, so that edge $j$ of $Q$
appears next to the left hand end of string $j$, for each $j$. We
write $Q$ down the right hand side also.

\begin{figure}
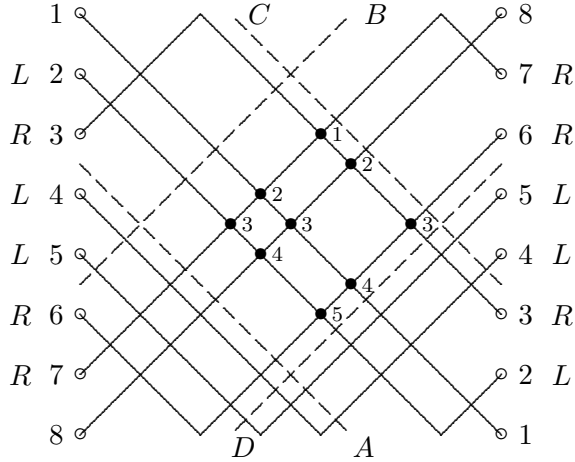

\beginpicture

\setcoordinatesystem units <0.8cm,0.8cm>             
\setplotarea x from -7 to 9, y from -.5 to 8       

\linethickness=0.5pt           

\multiput{$\circ$} at 0 0 *7 0 1 /
\multiput{$\circ$} at 7 0 *7 0 1 /

\put{$8$} at -0.4 0 %
\put{$7$} at -0.4 1 %
\put{$6$} at -0.4 2 %
\put{$5$} at -0.4 3 %
\put{$4$} at -0.4 4 %
\put{$3$} at -0.4 5 %
\put{$2$} at -0.4 6 %
\put{$1$} at -0.4 7 %

\put{$1$} at 7.4 0 %
\put{$2$} at 7.4 1 %
\put{$3$} at 7.4 2 %
\put{$4$} at 7.4 3 %
\put{$5$} at 7.4 4 %
\put{$6$} at 7.4 5 %
\put{$7$} at 7.4 6 %
\put{$8$} at 7.4 7 %

\put{$R$} at -1 1 %
\put{$R$} at -1 2 %
\put{$L$} at -1 3 %
\put{$L$} at -1 4 %
\put{$R$} at -1 5 %
\put{$L$} at -1 6 %

\put{$L$} at 8 1 %
\put{$R$} at 8 2 %
\put{$L$} at 8 3 %
\put{$L$} at 8 4 %
\put{$R$} at 8 5 %
\put{$R$} at 8 6 %

\setlinear \plot 0 7  7 0 / 
\setlinear \plot 0 6  6 0 / %
\setlinear \plot 6 0  7 1 / 
\setlinear \plot 0 5  2 7 / %
\setlinear \plot 2 7  7 2 / 
\setlinear \plot 0 4  4 0 / %
\setlinear \plot 4 0  7 3 / 
\setlinear \plot 0 3  3 0 / %
\setlinear \plot 3 0  7 4 / 
\setlinear \plot 0 2  2 0 / %
\setlinear \plot 2 0  7 5 / 
\setlinear \plot 0 1  6 7 / %
\setlinear \plot 6 7  7 6 / 
\setlinear \plot 0 0  7 7 / 

\setdashes <2mm,1mm>
\setlinear \plot 0 4.5  4.5 0 / %
\setlinear \plot 0 2.5  4.5 7 / %
\setlinear \plot 7 4.5  2.5 0 / %
\setlinear \plot 7 2.5  2.5 7 / %

\put{$\bullet$} at 2.5 3.5 %
\put{$\bullet$} at 3 4 %
\put{$\bullet$} at 4 5 %

\put{$\bullet$} at 3 3 %
\put{$\bullet$} at 3.5 3.5 %
\put{$\bullet$} at 4.5 4.5 %

\put{$\bullet$} at 4 2 %
\put{$\bullet$} at 4.5 2.5 %
\put{$\bullet$} at 5.5 3.5 %

\scriptsize{
\put{$3$} at 2.65 3.5 %
\put{$2$} at 3.15 4 %
\put{$1$} at 4.15 5 %

\put{$4$} at 3.15 3 %
\put{$3$} at 3.65 3.5 %
\put{$2$} at 4.65 4.5 %

\put{$5$} at 4.15 2 %
\put{$4$} at 4.65 2.5 %
\put{$3$} at 5.65 3.5 %
}

\normalsize{
\put{$A$} at 4.3 -0.2 %
\put{$B$} at 4.5 7 %
\put{$C$} at 2.57 7 %
\put{$D$} at 2.3 -0.2 %

}
\endpicture
\caption{Example: $Q=RRLLRL, Y=--LL--$.\label{rho}}
\end{figure}

Draw $4$ dashed lines on the diagram. The first, $A$,
should intersect the left vertical boundary between strings $a$ and $a+1$,
and go down and to the right at an angle of $-\pi/4$ (to the horizontal).
The next, $B$,
should intersect the left vertical boundary between strings $b$ and $b-1$
and go up and to the right at an angle of $\pi/4$. Line $C$ should
intersect the right vertical boundary between strings $a$ and $a+1$ and
go up and to the left at an angle of $3\pi/4$. Finally, line $D$
should intersect the right vertical boundary between strings $b$ and $b-1$,
and go down and to the left, at an angle of $-3\pi/4$.
For an example, in type $A_7$, where $Q=RRLLRL$ and $Y=--LL--$, with
$a(Y)=3$ and $b(Y)=6$, see Figure~\ref{rho}.

These four dashed lines bound a rectangle in the diagram; we shall see that
the crossings inside it form a rectangle, $\rho$, and that their labels are
precisely the positive roots in $M(Y)$. (In the example, these
crossings are marked by filled circles.) 

By construction, the lines crossing $\rho$ parallel to $A$ and $C$
are parts of the strings $c$ with $c\leq a$, and the
lines crossing $\rho$ parallel to $B$ and $D$ are parts of the strings
$d$ with $b\leq d$. Thus the crossings appearing inside $\rho$ are those
labelled with roots $\alpha_{cd}$ with $1\leq c\leq a\leq b\leq d\leq n+1$,
that is, the roots in $M(Y)$.

Let $E$ be the line in the rectangle we are interested in which is parallel
to $C$ and closest to $C$, and suppose it is part of string $x$ and
that this corresponds to an $R$ in $Q$. Let $y$
be the first string that $x$ intersects with after changing direction at
the top of the diagram. By the construction of the diagram, the strings
corresponding to the $L$'s of $Y$ are below the rectangle, so $y$ cannot
be one of these strings. Similarly, if $y$ were a string corresponding to an
$L$ to the right of $Y$, it would be parallel to $E$ initially (i.e. in the
left part of the diagram) and then would only be perpendicular to $E$ below
$\rho$. Also, if $y$ corresponded to an $R$ to the right of $Y$ in $Q$, $x$
and $y$ would intersect to the left of the point where $x$ changed
direction. We conclude that $x$ and $y$ intersect inside $\rho$.

Thus in the corresponding chamber diagram, the top crossing of the
rectangle is on row $1$ (so the corresponding simple reflection is $s_1$).
Suppose next that $x$ corresponds to an $L$ in $Q$. Let $y$ be the first
string that $x$ intersects. Since there can in this case be no $R$'s to
the right of $Y$ in $Q$, $x$ intersects every string corresponding to an
arrow of $Q$ to the right of $Y$ below the rectangle $\rho$; similarly with
the $L$'s in $Y$. Thus $y$ is to the left of $P$ and $x$ and $y$ intersect
in the rectangle.

Thus in both cases the top corner of the rectangle $\rho$ is in row $1$.
Let us
write the row of each crossing in $\rho$ next to it. The number of crossings
in a side of the rectangle parallel to $A$ or $C$ (equal to the number
of crossings involving $E$ inside $\rho$) is equal to $n+2-b$, and the
number in a side parallel to $B$ or $D$ is $a$. Thus the labels of the
crossings inside $\rho$ are the same and in the same configuration as the
numbers in $\rho(Y)$.

The case when $Y$ is of type $R$ is similar; in this
case the bottom corner of the rectangle in the arrangement diagram is
labelled $n$.

Now suppose that $Y$ is of type $L$ and immediately to the right of $Y$ in
$P$ is a component $Z$ of type $R$. By the construction of
the rectangles, the leftmost crossings inside them coincide. Similarly, if
$Y$ is of type $R$ and $Z$ is of type $L$, the rightmost crossings in the
rectangles coincide.

Thus the structure of the rectangles corresponding to components
of $P$ in the arrangement diagram, or the chamber diagram,
corresponding to $\bf i$, is the same as in $D(P)$. The result follows
by Theorem~\ref{multiset}.

In the case where $P=j\in [1,n]$, let $Y=j$; here
we consider $Y$ to be the vertex $j$ of $Q$, with no edges.
We then argue as in the above proof for $Y$ of type $L$.~$\square$

\begin{cor} \label{spanningmonomialcor}
Suppose we are in the situation of the Theorem, and let $\widetilde{F}(P)$ be
the monomial in the root operators as defined in~\ref{diagonalrow}.
Then $\widetilde{F}_{i_1}^{a_1}\widetilde{F}_{i_2}^{a_2}\cdots \widetilde{F}_{i_k}^{a_k}
\cdot 1=\widetilde{F}(P)$.
\end{cor}

\noindent {\bf Proof}: This follows immediately from the Theorem, since
the only relations in $U^-$ that are used are commutations, which the
root operators $\widetilde{F}_i$ also satisfy.~$\square$

\section{A Reparametrization}
We will now use the rectangle combinatorics from the previous section 
to calculate the Lusztig parametrizations
$c_j(\mathbf{j})$, $c_P(\mathbf{j})$ of the canonical basis elements
arising in Theorem~\ref{natural}; here $\mathbf{j}$ is fixed as the
reduced expression $(n,n-1,n,n-2,n-1,n,\ldots ,1,2,\ldots ,n)$ for $w_0$.


\begin{definition} \label{centralline} \rm
Suppose $P$ is a partial quiver and
$D(P)$ is the diagram for $P$ as defined above.
The rectangles
$\rho(Y)$ for $Y$ a component of $P$ divide $D(P)$ into
smaller rectangles, which we shall call boxes. (In the example near the
start of Section~\ref{rectangles}, there are $6$ boxes).
We regard these boxes as appearing in diagonal
rows, from top left to bottom right
As we work our way down the diagonal
rows of boxes, starting from the top right of $D(P)$,
the number of boxes in each diagonal is
at first odd, and then even, or is at first even, and then odd. We draw a
central line $Z$ on $D(P)$ dividing these two sets of diagonal rows of
boxes. In our example, the number of boxes in each row is $1$,$3$,$2$, so
the central line is indicated by the dashed line in Figure~\ref{centrallinef}.

\begin{figure}
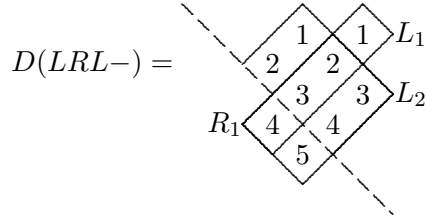

\beginpicture

\setcoordinatesystem units <0.8cm,0.8cm>             
\setplotarea x from -9 to 3, y from -2 to 3       

\linethickness=0.5pt           


\put{$D(LRL-)=$} at -2.5 1 %

\put{$4$}[c] at 0.5 0 %
\put{$3$}[c] at 1 0.5 %
\put{$2$}[c] at 1.5 1.0 %
\put{$1$}[c] at 2 1.5 %

\setlinear \plot 0 0  2 2 / %
\setlinear \plot 2 2  2.5 1.5 / %
\setlinear \plot 2.5 1.5  0.5 -0.5 / %
\setlinear \plot 0.5 -0.5  0 0 / %


\put{$5$}[c] at 1 -0.5 %
\put{$4$}[c] at 1.5 0 %
\put{$3$}[c] at 2 0.5 %

\setlinear \plot 0 0  1.5 1.5 / %
\setlinear \plot 1.5 1.5  2.5 0.5 / %
\setlinear \plot 2.5 0.5  1 -1 / %
\setlinear \plot 1 -1  0 0 / %


\put{$2$}[c] at 0.5 1 %
\put{$1$}[c] at 1 1.5 %

\setlinear \plot 0 1  1 2 / %
\setlinear \plot 1 2  2.5 0.5 / %
\setlinear \plot 2.5 0.5  1.5 -0.5 / %

\setdashes <2mm,1mm>
\setlinear \plot -1 2 2.5 -1.5 / %

\put{$L_1$}[c] at 2.8 1.5 %
\put{$R_1$}[c] at -0.3 0 %
\put{$L_2$}[c] at 2.8 0.5 %


\endpicture
\caption{The central line in $D(LRL-)$.\label{centrallinef}}
\end{figure}

Let $Y$ be a component of $P$, and define a set $S(Y)$ of
positive roots as follows. The rectangle $\rho(Y)$ for $Y$ is divided into
two parts by the central line. Let $\rho^*(Y)$ be the part of $\rho(Y)$
which is labelled.
The root $\alpha_{ij}$ lies in $S(Y)$ if and only
if the intersection of a diagonal row of $D(P)$ (see
Definition~\ref{diagonalrow}) with $\rho^*(Y)$ contains the numbers
$i,i+1,\ldots ,j-1,$ in order, reading from top left to bottom right.
Then $S(P)$ is defined to be the union of the sets $S(Y)$ for $Y$ a
component of $P$.
In our example, $S(L_1)=\{\alpha_1,\alpha_2,\alpha_3\}$,
$S(R_1)=\{\alpha_4+\alpha_5\}$ and $S(L_2)=\{\alpha_1+\alpha_2+\alpha_3,
\alpha_2+\alpha_3+\alpha_4\}$, so $S(P)=\{\alpha_1,\alpha_2,\alpha_3,
\alpha_4+\alpha_5,\alpha_1+\alpha_2+\alpha_3,\alpha_2+\alpha_3+\alpha_4\}$.

Let $v(P)\in \mathbb{N}^k$ be defined as follows. The reduced expression
$\bf j$ defines an ordering on the positive roots; we write $\beta^l=
s_{j_1}s_{j_2}\cdots s_{j_{l-1}}(\alpha_{j_l})$, for $l=1,2,\ldots ,k$.
Then let $v(P)_l$ be equal to $1$ if $\beta^l\in S(P)$ and $0$ otherwise.
So in the example above, $v(P)=(0,1,0,0,0,1,0,1,0,1,0,0,1,0,1)$.

If $j\in [1,n]$, we define $S(j)$ by specifying that the root $\alpha_{ij}$
lies in $S(j)$ if and only if there is a diagonal row of
$D(j)$ containing the numbers $i,i+1,\ldots ,j-1,$ in order, reading from
top left to bottom right. We define $v(j)$ in the same way as $v(P)$
is defined using $S(P)$.
\end{definition}

We know that, if $b\in {\bf B}$ and $i\in [1,n]$, then $\widetilde{F}_ib\equiv
b'\ ({\rm mod\ } vL')$ for some unique $b'\in {\bf B}$. With respect to
our fixed
reduced expression ${\bf j}$, this defines an action on $\mathbb{N}^k$,
defined by $\widetilde{F}_i(\varphi_{\bf j}(b))=\varphi_{\bf j}(b')$.
Using the identification in the last paragraph, this also induces an action
on sets of positive roots, which we also denote by $\widetilde{F}_i$.
Reineke, in~\cite[\S1]{reineke1},
describes the action of $\widetilde{F}_i$ on $\mathbb{N}^k$:

\begin{prop} (Reineke) \label{reineke} \\
Suppose ${\bf v}=(v_1,v_2,\ldots ,v_k)\in \mathbb{N}^k$. Write, for $1\leq i<
j\leq n+1$, $v_{ij}=v_l$, where $\beta^l=\alpha_{ij}$.
For such $i<j$, define
$$f_{ij}=\sum_{l=j}^{n+1} v_{il}-\sum_{l=j+1}^{n+1} v_{i+1,l}.$$
Let $j_0$ be minimal so that $f_{ij_0}=\max_j f_{ij}$. Then $\widetilde{F}_i$
increases $v_{ij_0}$ by $1$, decreases $v_{i+1,j_0}$ by $1$ (unless
$i+1=j_0$, when this latter effect does not occur), and leaves the other
$v_{ij}$'s unchanged.
\end{prop}

{\bf Definition:} Let $2\leq a\leq b\leq n$, and let $Q$ be a quiver.
We define the {\em $(a,b)$-sub partial quiver} of $Q$ to be the partial
quiver $P$ with $P\leq Q$, leftmost directed edge numbered $b$ and
rightmost directed edge numbered $a$.

As conjectured by Carter, we have:

\begin{theorem} \label{Scorrespondence}
Let $Q$ be a quiver, with ${\bf i}\in \chi_n$ compatible with $Q$. Suppose
that $(c,{\bf i})$ is a chamber of $Q$, with $P=P(c,{\bf i})$.
Let ${\bf j}=(n,n-1,n,n-2,n-1,n,\ldots ,1,2,\ldots ,n)$.
Let $b_P$ be the canonical basis element arising in Theorem~\ref{natural}.
Then the Lusztig parametrization $c_P(\mathbf{j})$ of $b_P$ is given by
$v(P)$. Similarly the Lusztig parametrization $c_j(\mathbf{j})$ of $b_j$ is
given by $v(j)$.
\end{theorem}

\noindent {\bf Proof:}
We have to calculate $S_{\bf i}^{\bf j}(a(c,{\bf i}))$ or
$S_{\bf i}^{\bf j}(a(j,{\bf i}))$, and thus we consider the corresponding
monomial of root operators acting on $1$,
$\widetilde{F}(P)=\widetilde{F}_{i_1}^{a_1}
\widetilde{F}_{i_2}^{a_2}\cdots \widetilde{F}_{i_k}^{a_k}\cdot 1$. Here
$P=P(c,{\bf i})$ in the first case and $j$ in the second.
Theorem~\ref{spanningmonomial} gives a description of
$\widetilde{F}(P)$ in terms of $D(P)$.
Recall that $D(P)$ consists of a number of {\em diagonal rows}
(see Definition~\ref{diagonalrow}), where diagonal row number $1$ is the
one which is at the top right of $D(P)$, diagonal row number $2$ is below
this (down and to the left), and so on. We now make a claim.

Let $x$ be a position in $D(P)$ containing a number, $i$. Suppose $x$ is
in (diagonal) row $a$. Let $r$ be the total number of rows in $D(P)$,
and let $\widetilde{F}(s)$, $s=1,2,\ldots ,r$, be the part of $\widetilde{F}(P)$
arising from the $s$th row. Thus $\widetilde{F}(P)=\widetilde{F}(r)\widetilde{F}(r-1)
\cdots \widetilde{F}(1)\cdot 1$.
For $s=1,2,\ldots ,r$, let $d_s$ be the number at
the rightmost end of the $s$th row. Define
$\widetilde{F}[x]:=\widetilde{F}_{i}^{\lambda_{i}}
\widetilde{F}_{i+1}^{\lambda_{i+1}}\cdots \widetilde{F}_{d_a}^{\lambda_{d_a}}
\widetilde{F}(a-1)\cdots \widetilde{F}(2)\widetilde{F}(1)\cdot 1$,
where, if $j$ appears in row $a$, $\lambda_j$ is the multiplicity
associated to it.

So $\widetilde{F}[x]$ is a final segment of $\widetilde{F}(P)$, arising from the
first $a-1$ rows of the diagram, together with the part of the $a$th row
to the right of $x$, including $x$.
Let ${\bf c}\in \mathbb{N}^k$ be defined by
$\widetilde{F}[x]\equiv F_{\bf i}^{\bf c}\ ({\rm mod\ } v{\cal L'})$.
Note that there is no problem working here modulo $v{\cal L'}$. We shall
have
$\widetilde{F}[x]\equiv b\ ({\rm mod\ } v{\cal L'})$ for a unique $b\in {\bf B}$
from
Kashiwara's approach, and $b\equiv F_{\bf i}^{\bf c}\ 
({\rm mod\ } v{\cal L'})$ for
a unique ${\bf c} \in \mathbb{N}^k$ by Lusztig's approach
and the comparison theorem~\cite[\S2.3]{lusztig3}.

Let $S(x)$ be the corresponding multiset of positive roots, obtained by
using the ordering on the set of positive roots induced by $\bf j$.
We claim that $S(x)$ is given by the union of the following two sets:
The first, $S_1(x)$ is defined to be those roots in $S(P)$ arising from rows
$1,2,\ldots ,a-1$, together with those roots in $S(P)$ arising from row $a$
of the form $\alpha_{s}+\alpha_{s+1}+\cdots +\alpha_t$ where $s\geq i$.
The second, $S_2(x)$, is defined to be the set of roots of form
$\alpha_{i}+\alpha_{i+1}+\cdots +\alpha_t$, where $\alpha_s+\alpha_{s+1}+
\cdots +\alpha_t$ is a root in $S(P)$ arising from row $a$ with $s<i$
(we call this procedure of passing from $\alpha_s+\alpha_{s+1}+\cdots +
\alpha_t$ to $\alpha_{i}+\alpha_{i+1}+\cdots +\alpha_t$, `cutting off').

We will prove the claim by induction on $x$. Let us assume that
for a given $x$ in $D(P)$, $S(x)$ is as described. Let $\bf c$ be as above,
corresponding to $x$. Let $y$ be one position
to the left of $x$ in the $a$th row of $D(P)$. We understand this to
mean the rightmost element of the $(a+1)$st row of $D(P)$ if $x$ is the
leftmost element of the $a$th row. In this case for the purposes of the
argument below we actually take $x$ to be
in the same row as $y$, immediately down and to the right of $y$, with
corresponding (fictitious) label $1$ more than that of $y$
(i.e. actually just outside the diagram $D(P)$).
Note that in this case, as we are assuming the claim to be true for $x$,
we have $S_2(x)=\phi$. Suppose there is an $i$
in position $y$. Then there is an $i+1$ in position $x$. Suppose that
$x$ and $y$ both lie in the $a$th row of $D(P)$.
We will show that the claim is also true
for $y$. The base case, when $x$ is taken to be
one position to the right of (and below) the rightmost position of the first
row, with $S(x)=\phi$, is trivially true. The Theorem will follow from the
induction argument, as it is the claim for
the case when $x$ is the leftmost position of the last row of the diagram.

We must show that the set of positive roots (defined by $\bf j$)
corresponding to $\bf c'$, where $\widetilde{F}[y]=\widetilde{F}_i^{\lambda_i}
\widetilde{F}[x]\equiv F_{\bf i}^{\bf c'}\ ({\rm mod\ } v{\cal L'})$ is given by
$S_1(y)\cup S_2(y)$, where $i$ is the number labelling position
$y$. In the case where $y$ is the rightmost element of its row (numbered $a$),
we take $\widetilde{F}[x]$ to be
$\widetilde{F}(a-1)\cdots \widetilde{F}(2)\widetilde{F}(1)\cdot 1$.
We define ${\bf c}\in \mathbb{N}^k$ by writing
$c_{ij}=1$ if $\alpha_{ij}\in S(x)$ and $c_{ij}=0$ otherwise. We
start with an important

\begin{lemma} \label{bothequalone}
If $l>i+1$ then we cannot have both $c_{il}=1$ and $c_{i+1,l}=1$.
\end{lemma}

\begin{figure}
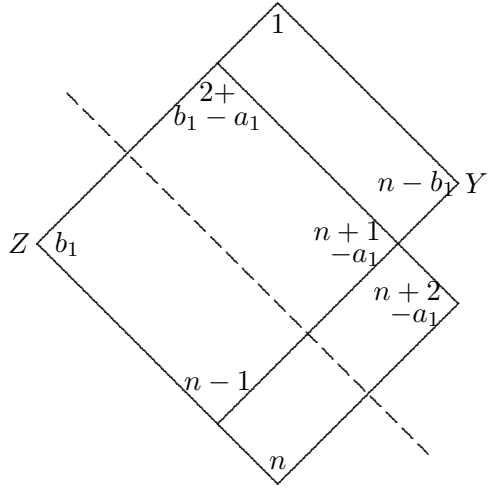

\beginpicture

\setcoordinatesystem units <0.8cm,0.8cm>             
\setplotarea x from -6.4 to 5, y from -4.5 to 5       

\linethickness=0.5pt           

\setlinear \plot 0 0 4 4 / %
\setlinear \plot 0 0 4 -4 / %
\setlinear \plot 3 3 7 -1 / %
\setlinear \plot 7 -1 4 -4 / %
\setlinear \plot 4 4 7 1 / %
\setlinear \plot 7 1 3 -3 / %
\setdashes <2mm,1mm>
\setlinear \plot 0.5 2.5 6.5 -3.5 / %

\put{$b_1$} at 0.5 0 %
\put{$2+$} at 3 2.5 %
\put{$b_1-a_1$} at 3 2.1 %
\put{$1$} at 4 3.65 %
\put{$n-b_1$} at 6.3 1 %
\put{$n+1$} at 5.15 0.2 %
\put{$-a_1$} at 5.28 -0.2 %
\put{$n-1$} at 3 -2.3 %
\put{$n+2$} at 6.15 -0.8 %
\put{$-a_1$} at 6.28 -1.2 %
\put{$n$} at 4 -3.65 %

\put{$Z$}[c] at -0.3 0 %
\put{$Y$}[c] at 7.3 1 %


\endpicture
\caption{The rectangles $\rho(Y)$ and $\rho(Z)$.\label{rhoyrhoz}}
\end{figure}

\noindent {\bf Proof:} (of Lemma).
First of all note that if $P=j\in [1,n]$, the Lemma is obvious, following
easily from the structure of $D(P)$. So we can assume that $P$ is a partial
quiver.
Suppose, for a contradiction, that
$c_{il}=c_{i+1,l}=1$. So $\beta_1:=\alpha_i+\alpha_{i+1}+\cdots +
\alpha_{l-1}$ and $\beta_2:=\alpha_{i+1}+\alpha_{i+2}+\cdots +\alpha_{l-1}$
both lie in
$S(x)$. We are assuming the claim to be true for $x$, so $S(x)=S_1(x)\cup
S_2(x)$. We shall see that the structure of $D(P)$ leads to a contradiction;
to do this we consider various cases: \\
Case (a) For $i=1,2$, if $\beta_i\in S_2(x)$, let
$\gamma_i$ be the root in $S(P)$ which is `cut off' to give $\beta_i$.
Otherwise, let $\gamma_i=\beta_i$.
Suppose that $\gamma_1,\gamma_2$ both come a rectangle for a component
$Y$ of $P$ of type $L$ which is the $(a_1,b_1)$-sub partial
quiver of $Q$. The set of positive roots in $S(P)$ arising from $Y$ is
$\{\alpha_{1,n-b_1},\alpha_{2,n+1-b_1},\ldots ,\alpha_{s,n+s-1+b_1}\}$
(where the central line passes through the rectangle corresponding to $Y$
immediately below the row with leftmost number $s$).
If $\alpha=\alpha_{pq}$ is any positive root, write $r(\alpha)=q-1$ and
$l(\alpha)=p$.
For each $\gamma$ in this list, $r(\gamma)$ is distinct, and since `cutting
off' a root $\gamma$ does not affect $r(\gamma)$, we conclude that this
case cannot occur, since $r(\beta_1)=r(\beta_2)$. \\
Case (b) Suppose that $\beta_1,\beta_2$ both come from a component
of $P$ of type $R$, as in (a). We argue as in case (a) to see that
this case cannot occur. \\
Case (c) Suppose that $\beta_1$ and $\beta_2$ arise from distinct components
of $P$ and that $\beta_1,\beta_2\in S_1(x)$.
Let $h$ be the height of the roots in $S(P)$ arising from a component
$Y$ of $P$. Then the height of the roots in $S(P)$ arising
from the component of $P$ immediately to the right of $Y$ (if
such exists) is $h+t$, where $t$ is the number of directed edges in $Y$.
Thus, since $\mbox{ht}(\beta_1)=1+\mbox{ht}(\beta_2)$ (ht denotes height),
for this case to occur, we must be in the situation where we have two
components $Y$ and $Z$ of $P$, where $Z$ immediately
follows $Y$, and $Y$ has precisely one directed edge. We assume that
$Y$ is of type $L$ and $Z$ is of type $R$ (the argument in the other case
is very similar). Suppose that $Z$ is the $(a_1,b_1)$-sub partial quiver of
$Q$. The two rectangles $\rho(Y)$ and $\rho(Z)$ will appear in
$D(P)$ as in Figure~\ref{rhoyrhoz},
where as usual, the dashed line indicates the central line in $D(P)$. Then
the roots for $Y$ will be taken from above this line, and the roots for
$Z$ from below, such roots $\alpha$ will have differing values of $r(\alpha)$
so we see that this case cannot occur. \\
Case (d) Suppose we have $\beta_1,\beta_2\in S_2(x)$. Suppose first that $x$
lies above the central line. We know that the lines in $D(P)$ which are the
bottom right edges of rectangles $\rho(Y)$ where $Y$ is a component
of $P$ of type $L$ are all distinct and parallel. It follows that
$r(\beta_1)\not=r(\beta_2)$, a contradiction. A similar argument holds if
$x$ lies below the central line. \\
Case (e) Suppose that $\beta_1\in S_1(x)$ and $\beta_2\in S_2(x)$. If
$\beta_1$ arises from the $a$th row of $D(P)$, we can argue as in (d)
(we will have again $r(\beta_1)\not= r(\beta_2)$, a contradiction).
So assume that $\beta_1$ arises from the $b$th row, where $1\leq b<a$.
We know that $\beta_2$ arises via cutting off from a root $\gamma_2\in S(P)$
which comes from the $a$th row of $D(P)$. Since $\mbox{ht}(\beta_2)<
\mbox{ht}(\gamma_2)$ and $\mbox{ht}(\beta_1)=1+\mbox{ht}(\beta_2)$, we
have $\mbox{ht}(\beta_1)< \mbox{ht}(\gamma_2)$. Therefore, arguing
as in (c), we see that $\beta_1$ comes from a component
$X_1$ occurring to the left of that giving rise to $\gamma_2$
(and thus $\beta_2$). We call the latter $X_2$. Suppose that $X_1$ does not
occur immediately to the left of $X_2$, or that $X_2$ is of type $R$.
Then by the
structure of $D(P)$, the bottom right lines $T_1$, $T_2$ delimiting the
rectangles $\rho(X_1)$ and $\rho(X_2)$ are distinct. Indeed, the line for
$\rho(X_2)$ occurs below and to the right of $\rho(X_1)$.

We know that there is an $l-1$ in the intersection of the $b$th row with
$\rho(X_1)$,
immediately to the left and above the line $T_1$. Similarly, there is an
$l-1$ in the intersection of the $a$th row with $\rho(X_2)$, immediately to
the left and above the line $T_2$, by the construction of $\beta_1,\beta_2$.
The relative position of lines $T_1$ and $T_2$ tells us that the former
$l-1$ lies strictly to the left of the latter, while the fact that $b<a$
tells us that the latter $l-1$ lies strictly to the left of the former, a
contradiction. We are left with the case when $X_1$ immediately precedes
$X_2$ and $X_2$ is of type $L$. In this case $T_1$ and $T_2$ coincide along
part of their lengths. Now $X_1$ is of type $R$, so $\beta_1$ arises from
a row below the central line, and $X_2$ is of type $L$, so $\beta_2$
arises from a row above the central line. It follows that
$r(\beta_1)\not=r(\beta_2)$, a contradiction. \\
Case (f) Suppose that $\beta_1\in S_2(x)$, $\beta_2\in S_1(x)$, and $\beta_1$
and $\beta_2$ arise from different components. If
$\beta_2$ arises from the $a$th row, we can argue as in (d). So assume
$\beta$ arises from the $b$th row, where $1\leq b<a$. Then $\beta_1$
arises via cutting off from a root $\gamma_1\in S(P)$ which comes from the
$a$th row. Since $\mbox{ht}(\beta_1)<\mbox{ht}(\gamma_1)$ and
$\mbox{ht}(\beta_1)=1+\mbox{ht}(\beta_2)$, we have $\mbox{ht}(\beta_2)<
\mbox{ht}(\gamma_1)$, and we can argue as in (e).

We have now covered all possible cases for $\beta_1,\beta_2\in S(x)$, and
we see that the Lemma is proved.~$\square$

We now investigate further the structure of $\bf c$. By assumption, for
each $j,k$, $c_{jk}$ is equal to $0$ or $1$. The above shows that we
cannot have both $c_{il}=c_{i+1,l}=1$ for $l>i+1$. We shall now pin down
further the possible values for $c_{il}$ and $c_{i+1,l}$. We know that,
for each $l>i+1$, $(c_{il},c_{i+1,l})=(0,0)$, $(0,1)$ or $(1,0)$. If we
let $c_{i+1,i+1}:=1-c_{i,i+1},$ then this is also true for $l=i+1$.
For $m=0,1,2,\ldots n-i$,
write $a_m=c_{i,i+m+1}$ and $b_m=c_{i+1,i+m+1}$, so
the pairs we are interested in are of form ${\bf x}_m=(a_m,b_m)$,
$m=0,1,2,\ldots n-i$. We first need a Lemma describing the action of
$\widetilde{F}_i$ in terms of these pairs, using Proposition~\ref{reineke}:

\begin{lemma} \label{applyreineke}
Suppose $\bf c$ is as in the above, and we have written the pairs
${\bf x}_0,{\bf x}_1,\ldots ,{\bf x}_{n-i}$ as above, with each ${\bf x}_m$
equal to $(0,0)$, $(0,1)$ or $(1,0)$. Rewrite this sequence, replacing
pairs of these types with $0$, `$-$' and `$+$', respectively. Cross out all
$0$'s, except initial $0$'s. Cross out all non-initial pairs of symbols of
form $+-$, and repeat until there are no such pairs of symbols in the
sequence left. Use the following table of possible sequences to determine
a symbol in the resulting sequence of symbols. This is defined
to be the symbol that is indicated by a line above it. Every time a
sequence of symbols appears in the list, as in $-,-,\ldots ,-$,
this means that the symbol concerned must occur at least once.

\noindent (i) $0,-,-,\ldots ,\overline{-},+,+,\ldots ,+,$  \\
(ii) $0,-,-,\ldots ,\overline{-},$  \\
(iii) $\overline{0},+,+,\ldots ,+,$  \\
(iv) $\overline{0},$ \\
(v) $(+,-),-,-,\ldots ,\overline{-},+,+,\ldots ,+,$ \\
(vi) $(+,-),-,-,\ldots ,\overline{-},$ \\
(vii) $(\overline{+},-),+,+,\ldots ,+,$ \\
(viii) $\overline{+},-,$ \\
(ix) $-,-,\ldots ,\overline{-},+,+,\ldots ,+,$ \\
(x) $-,-,\ldots ,\overline{-},$ \\
(xi) $\overline{+},+,\ldots ,+$.

In this way, a symbol in the original sequence is determined, and therefore
an ${\bf x}_m$. Then acting $\widetilde{F}_i$ has the following effect on
$\bf c$. If ${\bf x}_m=(0,1)$ (i.e. a `$-$'), then $c_{i+1,i+m+1}=1$. This is
changed to zero (having no effect if $m=0$), and $c_{i,i+m+1}$ is changed
from zero to $1$, thus changing the corresponding symbol in the sequence
from a `$-$' to a `$+$'. If ${\bf x}_m=(1,0)$ or $(0,0)$, then we must have
$m=0$ (as no initial symbols are crossed out), and $\widetilde{F}_i$ has the
effect of increasing $c_{i,i+1}$ by $1$.
\end{lemma}

\noindent {\bf Proof}: We consider how the $f_{ij}$'s are related to the
${\bf x}_p$'s (see Proposition~\ref{reineke} for the definition of the
$f_{ij}$'s). 
For each $0\leq p\leq n-i$, write
$\delta_p:=f_{i,i+p+2}-f_{i,i+p+1}$.
A simple calculation shows us that
$\delta_p=1$ if $({\bf x}_p,{\bf x}_{p+1})=(0,-)$ or $(-,-)$, that
$\delta_p=-1$ if $({\bf x}_p,{\bf x}_{p+1})=(+,+)$ or $(+,0)$, and is
zero in all other cases. Thus the sequence in the Lemma enables us to
calculate the value of the $f_{ij}$'s, given the value of $f_{i,i+1}$ (for
our fixed $i$). Recall that Proposition~\ref{reineke} tells us that, in
order to see how $\widetilde{F}_i$ acts, we need to calculate the smallest $j$
for which $f_{ij}$ is maximal; write $j_0$ for this value of $j$.
It is easy to check that the crossing out routines described above never cross out the symbol corresponding to
${\bf x}_p$ where $i+p+1=j_0$. Thus recalculating this {\em after} doing
the crossing out does not affect the final result. After crossing out is
done, the only possible sequences are those as listed in the Lemma. It is
then easy to check that $j_0$ is given by $j_0=i+p+1$ where $p$ is given by
the rules in the Lemma in each case (i)-(xi).~$\square$

\noindent {\bf Example:} To illustrate this Lemma we go back to our example
(see Definition~\ref{diagonal} and Definition~\ref{centralline}). So
$P=LRL-$ and $D(P)$ is as in Figure~\ref{centrallinef}.
Suppose $x$ is the position of the
rightmost $4$ in $D(P)$, the last element in the third diagonal row of
$D(P)$. Then $y$ is the leftmost $3$ in $D(P)$, the second element in the
third diagonal row, so $i=3$. The set $S_1(x)$ is by assumption $\{\alpha_1,
\alpha_2,\alpha_1+\alpha_2+\alpha_3\}$, and $S_2(x)$ is $\{\alpha_4\}$.
Thus $S(x)$ is the union of these two sets and ${\bf c}=(0,0,1,0,0,0,0,0,0,
1,0,0,1,0,1)$. We can write $\bf c$ as in Figure~\ref{vectorc}.
\begin{figure}
\beginpicture

\setcoordinatesystem units <0.75cm,0.75cm>             
\setplotarea x from -4 to 5, y from 1 to 5.5       


\put{$c_{12}$} at 1 1 %
\put{$c_{23}$} at 2 1 %
\put{$c_{34}$} at 3 1 %
\put{$c_{45}$} at 4 1 %
\put{$c_{56}$} at 5 1 %
\put{$c_{13}$} at 1.5 2 %
\put{$c_{24}$} at 2.5 2 %
\put{$c_{35}$} at 3.5 2 %
\put{$c_{46}$} at 4.5 2 %
\put{$c_{14}$} at 2 3 %
\put{$c_{25}$} at 3 3 %
\put{$c_{36}$} at 4 3 %
\put{$c_{15}$} at 2.5 4 %
\put{$c_{26}$} at 3.5 4 %
\put{$c_{16}$} at 3 5 %

\put{$=$} at 7 3 %

\put{$1$} at 9 1 %
\put{$1$} at 10 1 %
\put{$\underline{0}$} at 11 1 %
\put{$\underline{1}$} at 12 1 %
\put{$0$} at 13 1 %
\put{$0$} at 9.5 2 %
\put{$0$} at 10.5 2 %
\put{$\underline{0}$} at 11.5 2 %
\put{$\underline{0}$} at 12.5 2 %
\put{$1$} at 10 3 %
\put{$0$} at 11 3 %
\put{$\underline{0}$} at 12 3 %
\put{$0$} at 10.5 4 %
\put{$0$} at 11.5 4 %
\put{$0$} at 11 5 %


\endpicture
\caption{The vector $\bf c$.\label{vectorc}}
\end{figure}

We have ${\bf x}_0=(c_{34},c_{44})=(0,1)$,
${\bf x}_1=(c_{35},c_{45})=(0,1)$, and
${\bf x}_2=(c_{36},c_{46})=(0,0)$ --- these are indicated by underlining
in Figure~\ref{vectorc} (note that $c_{44}$ does not appear in the diagram).
Thus the sequence of symbols as in Lemma~\ref{applyreineke} is
$(-,-,0)$. In the reduction step, we remove the final zero to get $(-,-)$
and see we are in case (x)
of the Lemma, and that the final `$-$' is marked. This corresponds to
${\bf x}_1$ and the Lemma tells us that applying $\widetilde{F}_i$ changes
$c_{35}$ from $0$ to $1$ and $c_{45}$ from $1$ to $0$.
As a more complicated example, consider the sequence of symbols
$(0,0,+,+,-,-,+,+,+,0,-,0,-,-,+,-,-,-,-,0,0)$. We first of all remove
non-initial zeros to get $(0,+,+,-,-,+,+,+,-,-,-,+,-,-,-,-)$ and then remove
$+,-$ pairs to get $(0,-,-,-)$ --- which is case (ii) of the Lemma. Thus
the last `$-$' is indicated, corresponding to the last `$-$' of the original
sequence.

We now begin the investigation as promised.
\begin{lemma} \label{investigation}
Suppose that ${\bf x}_j=(1,0)$. Then $c_{i,i+j+1}=1$, so
that $\alpha_{i,i+j+1}\in S(x)$. Suppose that in fact $\alpha_{i,i+j+1}\in
S_1(x)$. Then there is $m>j$ so that ${\bf x}_{j+1}={\bf x}_{j+2}=\cdots =
{\bf x}_{m-1}=(0,0)$ and ${\bf x}_m=(0,1)$. Thus $c_{i+1,i+m+1}=1$, so
$\alpha_{i+1,i+m+1}\in S(x)$; in fact $\alpha_{i+1,i+m+1}\in S_1(x)$.
\end{lemma}

\noindent {\bf Proof}:
Let $l=i+j+2$; we know that $\alpha_i+\alpha_{i+1}+
\cdots +\alpha_{l}\in S_1(x)$. So
$\alpha_i+\alpha_{i+1}+ \cdots +\alpha_{l}\in S(P)$ and arises from the
intersection of row $b$ of $D(P)$ with a rectangle $\rho$, where $b<a$ (since
$\alpha_i$ appears in it). Thus if $P=j\in [1,n]$, we have a contradiction
by the structure of $D(P)$ (which consists of only one rectangle); $i$
can appear at the start of one row in $D(P)$ only. So
in this case the circumstances of the Lemma cannot occur. Therefore we can
assume $P$ is a partial quiver and
we define $X$ to be the component of $P$ corresponding to $\rho$.
Suppose that row $b$ is not immediately above the central line in $D(P)$.
Then row $b+1$ of $D(P)$ intersects $\rho,$ the rectangle for $X$.
This is because $i$ occurs in row $a$, in position $y$,
where $a>b$. If row $b$ were the bottom row in the
rectangle for $X$, $i$ would appear in the left hand corner of this
rectangle, and therefore in no later row in $D(P)$ (by the structure of
$D(P)$), a contradiction.
Thus the intersection of row $b+1$ with the rectangle for $X$ gives a root
$\alpha_{i+1}+\alpha_{i+2}+\cdots +\alpha_{l+1}\in S(P)$. Since $b+1\leq a$
and the label at position $x$ in $D(P)$ is $i+1$, this root lies in $S_1(x)$
(i.e. is not cut off), as required (with $m=j+1$).

So we are left in the case where row $b$ does lie immediately above the
central line of $D(P)$. Then $X$ must be of type $L$.
Suppose first that $X$ is the
rightmost component of $P$. Then $i$ appears in the
leftmost column (where a column is perpendicular to what we are calling a
row) of $D(P)$, and also is immediately above the central line.
The central line coincides with the bottom left delimiting line
of the rectangle for $X$. Consider now the rectangle for the component
$Y$ immediately to the left of $X$. This has
bottom right delimiting line
parallel to that for $X$, but length (from top left to bottom
right) less than that for $X$. Thus
its upper left delimiting line is down and to the
right of that for $X$, and $i$ is in the left hand
corner of the rectangle for $X$. Then
all of the upper left delimiting lines for
components other than $X$ appear down and to the right of that for
$X$, so $i$ cannot appear below the central line.
But as $b<a$, $x$ is below the central line, and therefore
so is $y$, which is labelled $i$, and we have a contradiction.
Thus there is a component
immediately to the right of $X$ in $P$; we call this
$Y$. Let $\alpha_{i+1}+\alpha_{i+2}+\cdots +\alpha_t$ (some $t>l+1$) in
$S(P)$ be the root arising from the intersection of the $(b+1)$st row of
$D(P)$ with the rectangle for $Y$. Then as before $\alpha_{i+1}+
\alpha_{i+2}+\cdots +\alpha_t\in S_1(x)$. Put $m=t-i$, and we have
$(a_m,b_m)=(c_{i,i+m+1},c_{i+1,i+m+1})=(0,1)$, as required. We
now have to show that ${\bf x}_{j+1}={\bf x}_{j+2}=\cdots =
{\bf x}_{m-1}=(0,0)$. That is, we have to show that \\
(a) $\alpha_i+\alpha_{i+1}+\cdots +\alpha_u\not\in S(x)$ for $l+1\leq u<t$,
and \\
(b) $\alpha_{i+1}+\alpha_{i+2}+\cdots +\alpha_u\not\in S(x)$ for
$l+1\leq u<t$.

We start with (a), and take $l+1\leq u<t$. Since all roots
$\beta\in S_2(x)$ satisfy $l(\beta)=i+1$, we have
$\alpha_i+\alpha_{i+1}+\cdots +\alpha_{u}\not \in S_2(x)$.
Therefore, it is enough to show that
$\alpha_i+\alpha_{i+1}+\cdots +\alpha_{u}\not\in S(P)$. But
$\mbox{ht}(\alpha_i+\alpha_{i+1}+\cdots +\alpha_{u})$ is strictly between
the heights of the roots in $S(P)$ arising from $X$ and $Y$, and so by our
description of the heights of the roots in $S(P)$ in case (c) of the proof of
Lemma~\ref{bothequalone}, we see that 
$\alpha_i+\alpha_{i+1}+\cdots +\alpha_{u}\not\in S(P)$ and (a) holds.

For (b), we consider again $l+1\leq u<t$. If $u\not=l+1$, the argument
used for (a) shows that $\alpha_{i+1}+\alpha_{i+2}+\cdots +\alpha_{u}\not
\in S(P)$. But if $u=l+1$, then this root has the same height as roots in
$S(P)$ arising from $X$. So if it were in $S(P)$, it would have arisen
from the intersection of the rectangle for $X$ with row $b+1$ of $D(P)$.
But this row is below the central line of $D(P)$, whence
$\alpha_{i+1}+\alpha_{i+2}+\cdots +\alpha_{l+1}$ does not arise
in this way, and so does not lie in $S(P)$.
Next, suppose that we had
$\alpha_{i+1}+\alpha_{i+2}+\cdots +\alpha_{u}\in S_2(x)$. Then, since
$x$ must be below the central line, $u$ must be immediately to the upper
left of the bottom right line delimiting a rectangle corresponding to a
component $Z$ of $P$ of type $R$. Since
$\alpha_{i+1}+\alpha_{i+2}+\cdots +\alpha_{u}$ is assumed to be in
$S_2(x)$, it must be the `cut off' from a root $\gamma=\alpha_p+\alpha_{p+1}+
\cdots +\alpha_u\in S(P)$, $p\leq i$, where $\gamma$ arises from the
rectangle corresponding to $Z$. Since $\mbox{ht}(\gamma)>\mbox{ht}(
\alpha_{i+1}+\alpha_{i+2}+\cdots +\alpha_{l+1}),$ $Z$ must lie strictly to
the right of $X$ in $P$. We know that
$\alpha_{i+1}+\alpha_{i+2}+\cdots +\alpha_{t}$ corresponds to the first
line of the rectangle for $Y$ below the central line in $D(P)$, so 
$Z$ is strictly to the right
of $Y$, as $u<t$. Thus the bottom right delimiting line for $Z$ is
strictly down to the right of that for $Y$.
Thus this $u$ lies to the right of the $u$ corresponding to the $\alpha_u$
in $\alpha_{i+1}+\alpha_{i+2}+\cdots +\alpha_{t}$ (in the $(b+1)$st row
of $D(P)$), and this $u$ lies below the central line. This is impossible,
as all $u$'s in $D(P)$ to the right of the $u$ in the $(b+1)$st row
(which is immediately below the central line), lie {\em above} the central
line. Thus $\alpha_{i+1}+\alpha_{i+2}+\cdots +\alpha_{u}\not\in S_2(x)$
and so it is not in $S$ for $l+1\leq u<t$. So we see (b) also holds, and
the Lemma is proved.~$\square$

\noindent {\bf Remark}
Note that all roots in $S_2(x)$ have a corresponding ${\bf x}_j$ of
form $(0,1)$. All roots in $S_1(x)$ arising from the $a$th row of
$D(P)$ will be of form $\alpha_p+\alpha_{p+1}+\cdots +\alpha_l$, for
$l\geq p\geq i+1$, so only if $p=i+1$ will they contribute to an
${\bf x}_j$; in this case we shall get a corresponding ${\bf x}_j$ of form
$(0,1)$.

We can now apply Lemma~\ref{applyreineke} in our situation. By 
Lemma~\ref{applyreineke}, after the crossing out process, our sequence looks
like $((+,-),-,-,\ldots ,-),$ $((0),-,-,\ldots ,-)$ or $((-),-,\ldots ,-)$.
By Lemma~\ref{applyreineke}, applying $\widetilde{F}_i$ will change the last
minus to a plus, and repeated applications will change the final $\lambda_i$
minuses to pluses. By the proof of Lemma~\ref{investigation}, it is clear that
all of the $-$'s in the sequence, except for those in brackets, correspond
to precisely those ${\bf x}_p$'s which correspond to roots of form
$\alpha_{i+1}+\alpha_{i+2}+\cdots +\alpha_l\in S_2(x)$.
	Suppose first that there is no edge of a rectangle between $x$ and
$y$. By structure of $D(P)$, $\lambda_i$ is the number of roots in $S_2(x)$.
Applying $\widetilde{F}_i$ $\lambda_i$ times thus replaces each root
$\alpha_{i+1}+\alpha_{i+2}+\cdots +\alpha_l\in S_2(x)$
by $\alpha_i+\alpha_{i+1}+\cdots +\alpha_l$
(see Lemma~\ref{applyreineke}).
	Next, suppose that there is an edge which is the bottom right edge
of a rectangle between $x$ and $y$. By the structure of $D(P)$, $\lambda_i$
is $1$ more than the number of roots in $S_2(x)$. Applying $\widetilde{F}_i$
$\lambda_i-1$ times thus replaces each root
$\alpha_{i+1}+\alpha_{i+2}+\cdots +\alpha_l\in S_2(x)$
by $\alpha_i+\alpha_{i+1}+\cdots +\alpha_l$
(see Lemma~\ref{applyreineke}). We must have $c_{i,i+1}=0$, so that the
initial symbol in our sequence is a `$-$', and we must have a sequence of
form $(-),-,-,\ldots ,-$. From this we can conclude that applying
$\widetilde{F}_i$ one more time changes this initial `$-$' (in brackets)
into a `$+$', introducing the root $\alpha_i$ into $S(x)$.
	Finally, suppose there is an edge which is the top left edge of a
rectangle between $x$ and $y$. By the structure of $D(P)$, $\lambda_i$ is
$1$ less than the number of roots in $S_2(x)$. Applying $\widetilde{F}_i$
$\lambda_i$ times thus replaces each root
$\alpha_{i+1}+\alpha_{i+2}+\cdots +\alpha_l\in S_2(x)$
by $\alpha_i+\alpha_{i+1}+\cdots +\alpha_l$
(see Lemma~\ref{applyreineke}), except for the one with minimal value of $l$.
	Thus we see that in each case, applying $\widetilde{F}_i$ $\lambda_i$
times replaces $S(x)$ with $S(y)$, proving that our claim
is also true for $y$, and we see it is true for all points in the diagram.
It follows that Theorem~\ref{Scorrespondence} is true.~$\square$

\noindent {\bf Acknowledgements} \\
Thanks are due to Professor R. W. Carter of the University of Warwick,
whose help, ideas and discussions were invaluable.
The research for this paper was carried out while the author was an
EPSRC research assistant of Professor K.\,A. Brown at the University of
Glasgow, Scotland. Part of the work was completed while the author was
attending the programme `Representation Theory of Algebraic Groups and
Related Finite Groups' at the Isaac Newton Institute, Cambridge, in 1997.
The author is grateful to the referee for some helpful comments on an
earlier version of this manuscript. This paper was typeset in \LaTeX.

\newcommand{\noopsort}[1]{}\newcommand{\singleletter}[1]{#1}


\begin{thebibliography}{10}

\bibitem{bedard2}
R.~Bedard.
\newblock On the spanning vectors of {L}usztig cones.
\newblock Preprint, 1999.

\bibitem{bfz1}
A.~Berenstein, S.~Fomin, and A.~Zelevinsky.
\newblock Parametrizations of canonical bases and totally positive matrices.
\newblock {\em Adv. Math.}, {\bf 122}(1) (1996), 49--149, doi:10.1006/aima.1996.0057.

\bibitem{bz1}
A.~Berenstein and A.~Zelevinsky.
\newblock String bases for quantum groups of type {$A_r$}.
\newblock {\em Adv. Soviet Math.}, {\bf 16}, Part 1 (1993), 51--89.

\bibitem{bz3}
A.~Berenstein and A.~Zelevinsky.
\newblock Tensor product multiplicities, canonical bases and totally positive varieties.
\newblock Preprint {\tt arXiv:math.RT/9912012}, December 1999.

\bibitem{carter3}
R.~W.~Carter.
\newblock Canonical bases, reduced words, and {L}usztig's
piecewise-linear function.
\newblock {\em in} ``Algebraic Groups and Lie Groups'', editor G. I. Lehrer,
pp 61--79. C.U.P (1997).

\bibitem{me10}
R.~W. Carter and R.~J. Marsh.
\newblock Regions of linearity, {L}usztig cones and canonical basis elements
  for the quantized enveloping algebra of type ${A}_4$.
\newblock University of Leicester Technical Report 2000/17.

\bibitem{kash2}
M.~Kashiwara.
\newblock On crystal bases of the {$q$}-analogue of universal enveloping
  algebras.
\newblock {\em Duke Math. J.}, {\bf 63}(2) (1991), 465--516.

\bibitem{kash4}
M.~Kashiwara.
\newblock The crystal base and {L}ittelmann`s refined {D}emazure character
  formula.
\newblock {\em Duke Math. J.}, {\bf 71}(3) (1993), 839--858.

\bibitem{lusztig2}
G.~Lusztig.
\newblock Canonical bases arising from quantized enveloping algebras.
\newblock {\em J. Amer. Math. Soc.}, 3:447--498, 1990.

\bibitem{lusztig3}
G.~Lusztig.
\newblock Canonical bases arising from quantized enveloping algebras, {II}.
\newblock {\em Prog. Theor. Phys. Sul.}, {\bf 102} (1990), 175--201.

\bibitem{lusztig6}
G.~Lusztig.
\newblock ``Introduction to Quantum Groups,''
\newblock Birkh\"{a}user, Boston, 1993.

\bibitem{lusztig7}
G.~Lusztig.
\newblock Tight monomials in quantized enveloping algebras,
\newblock {\em in} ``Quantum deformations of algebras and their
representations,''
  Volume~7 of {\em Israel Math. Conf. Proc.}, pp. 117--132, 1993.

\bibitem{me8}
R.~J. Marsh.
\newblock The {L}usztig cones of quantized enveloping algebras of type {$A$}.
\newblock {Preprint {\tt arXiv:math.QA/0008035.}}

\bibitem{me7}
R.~J. Marsh.
\newblock More tight monomials in quantized enveloping algebras.
\newblock {\em J. Alg.}, {\bf 204} (1998), 711--732.

\bibitem{nz1}
T.~Nakashima and A.~Zelevinsky.
\newblock Polyhedral realizations of crystal bases for quantized {K}ac-{M}oody
  algebras.
\newblock {\em Adv. in Math.}, {\bf 131}(1) (1997), 253--278,
doi:10.1006/aima.1997.1670.

\bibitem{premat1}
A.~Premat.
\newblock The {L}usztig cone and the image of the {K}ashiwara map.
\newblock Private communication, 1999.

\bibitem{reineke1}
M.~Reineke.
\newblock On the coloured graph structure of {L}usztig's canonical basis.
\newblock {\em Math. Ann.}, {\bf 307}(4) (1997), 705--723.

\bibitem{reineke2}
M.~Reineke.
\newblock Monomials in canonical bases of quantum groups and quadratic forms.
\newblock Preprint, 1999.

\bibitem{xi3}
N.~Xi.
\newblock Canonical basis for type {$A_3$}.
\newblock {\em Comm. Alg.}, {\bf 27}(11) (1999), 5703--5710.

\bibitem{xi2}
N.~Xi.
\newblock Canonical basis for type {$B_2$}.
\newblock {\em J. Alg.}, {\bf 214}(1) (1999), 8--21.

\bibitem{xi4}
N.~Xi.
\newblock Private communication.
\newblock 1997.

\end{thebibliography}
\end{document}